\documentclass[10pt]{article}

\usepackage{latexsym}
\usepackage{amssymb}
\usepackage{amscd}
\usepackage{amsmath}
\usepackage{amsthm}
\usepackage[usenames,dvipsnames]{color}
\usepackage{graphicx}
\usepackage{wrapfig}
\usepackage{tikz}
\usetikzlibrary{shapes.geometric}
 \usepackage{url}
 \usepackage{fancyvrb}

\usepackage[mathscr]{euscript}
\usepackage{mathtools}

\usepackage{hyperref}
\hypersetup{
 colorlinks,
 citecolor=Violet,
 linkcolor=Blue,
 urlcolor=Blue}

\def\qed{\hfill\rule{2.0mm}{2.0mm}}
\def\pf{\noindent{\bf Proof:}~ }
\def\eop{\hfill\rule{2.0mm}{2.0mm}}

\def\bC{{{\mathbb C}}}
\def\bR{{{\mathbb R}}}

\def\geq{\geqslant}
\def\leq{\leqslant}

\newcommand{\goto}{\rightarrow}

\newcommand{\beq}{\begin{equation*}}
\newcommand{\eeq}{\end{equation*}}
\newcommand{\beqn}{\begin{equation}}
\newcommand{\eeqn}{\end{equation}}
\newcommand{\bea}{\begin{eqnarray}}
\newcommand{\eea}{\end{eqnarray}}

\newcommand{\gap}{\vspace{0.1in}}
\newcommand{\demo}[1]{\pf}
\newcommand{\bx}{{\bf x}}
\newcommand{\ba}{{\bf a}}

\newcommand{\Iso}{{\rm Iso}}

\newtheorem{proposition}{Proposition}[section]
\newtheorem{lemma}[proposition]{Lemma}
\newtheorem{theorem}[proposition]{Theorem}

\newtheorem{conjecture}[proposition]{Conjecture}
\theoremstyle{definition}

\theoremstyle{remark}
\newtheorem{remark}[proposition]{Remark}

\begin{document}
\title{\sc On the Uniqueness of Clifford Torus with Prescribed Isoperimetric Ratio}

\date{March 26, 2020}

\author{
 Thomas Yu\thanks{
Department of Mathematics, Drexel University. Email: \href{mailto:yut@drexel.edu}{yut@drexel.edu}.
He is supported in part by the National Science Foundation grants DMS 0512673 and DMS 0915068.}
\and
Jingmin Chen\thanks{
Citigroup Global Markets Inc., 390 Greenwich Street, New York, NY 10013, U.S.A.. Email: \href{mailto:jingmchen@gmail.com}{jingmchen@gmail.com}.
}
}

\makeatletter \@addtoreset{equation}{section} \makeatother
\maketitle

\centerline{\bf Abstract:}
The Marques-Neves theorem asserts that among all the torodial (i.e. genus 1) closed surfaces,
the Clifford torus has the minimal Willmore energy $\int H^2 \, dA$. 
Since the Willmore energy is invariant M{\"o}bius transformations, it can be shown that there is a one-parameter family, up to homotheties, of
genus 1 Willmore minimizers. It is then a natural conjecture that such a minimizer is unique if one prescribes its
isoperimetric ratio.
In this article, we show that this conjecture can be reduced to the positivity question of a polynomial recurrence.

\vspace{.2in} \noindent {\bf Acknowledgments.} This work is partially supported
by NSF grants DMS 0915068 and DMS 1115915. We thank Robert Kusner for bringing to
our attention the uniqueness problem. Also, we are grateful to
Manuel Kauers, Stephen Melczer and Pierre Lairez for sharing their expertise in P-recurrences.

\vspace{.2in}

\section{Uniqueness problem in the Canham-Evans-Helfrich model} \label{sec:intro}
Why do all humans of all races occur to have the same biconcave shaped red blood cells? This apparent uniqueness might
have intrigued biologists since the invention of microscope. The seminal work of Canham \cite{Canham:Elastic}, Helfrich \cite{Helfrich:Elastic} and
Evans \cite{Evans1974923} suggests that bending elasticity, induced by curvature, plays the key role in driving the
geometric configurations of such membranes.

The so-called spontaneous curvature model of Helfrich
suggests that a biomembrane surface $S$ configures itself to minimize
$\int_S H^2 dA$ subject to the area, volume and area difference (related to the bilayer characteristics) constraints,
i.e. $S$ solves the variational \textbf{Helfrich problem}
\begin{equation}
\begin{aligned}
 \min_S W(S) :=\int_S H^2 \, dA
\text{ s.t. }
\left\{
  \begin{array}{ll}
    \mbox{(i)} & A(S) := \int_S 1 \; dA = A_0, \\
    \mbox{(ii)} & V(S) :=
\frac{1}{3} \int_S [x {\hat{\mathbf{i}}} + y \hat{\mathbf{j}} + z \hat{\mathbf{k}}] \cdot \hat{\mathbf{n}} \: dA =
V_0, \\
    \mbox{(iii)} & M(S) :=-\int_S H \; dA = M_0.
  \end{array}
\right.
\end{aligned}
\label{eq:Helfrich}
\end{equation}
Here $H=(\kappa_1+\kappa_2)/2$ is the mean curvature. (We assume that the normal of any closed orientable surface points outward.
In particular, it means $H<0$ for a sphere.)
In (ii), $V(S)$ is the enclosed volume.
The connection of (iii)
to bilayer area difference comes from the
relation $-\int_S H dA = \lim_{\varepsilon \goto 0} \frac{1}{4\varepsilon}({\rm area} (S_{+\varepsilon})- {\rm area} (S_{-\varepsilon}))$,
where $S_{+\varepsilon}$ and $S_{-\varepsilon}$ are the `$\varepsilon$-offset surfaces',
 and
that the thickness of the lipid bilayer, $2 \varepsilon$, is negligible compared to the size of the vesicle.
The constraint values $A_0$, $V_0$ and $M_0$ are determined by physical conditions (e.g. temperature, concentration).
$W(S)$ is called the \textbf{Willmore energy} of the surface $S$.
When the area-difference constraint (iii) is omitted, the variational problem
 is referred to as the \textbf{Canham problem}.
When even the volume constraint (ii) is omitted, the variational problem is referred to as the \textbf{Willmore  problem}. In this case,
there is essential no constraint as $W$ is scale-invariant. In any case,
 the area constraint (i)
only fixes the scale; see the discussion around \eqref{eq:v0m0} below.

It is observed experimentally that no topological change occurs in any accessible time-scale, so
the Helfrich, Canham or Willmore problems ask for a minimizer $S$ over all orientable closed
surface with a fixed genus $g$. Spherical ($g=0$) vesicles are the most common among naturally
occurring biomembranes, although higher genus ones have been synthesized in the laboratory
\cite{MichaletBensimon:Genus2,LIPOWSKY:Genus2,Seifert:Config}.
The Canham, Helfrich and related models explain the large
variety of shapes observed in even a closed vesicle with a spherical topology.

At a mathematical level, the existence of solution for the Canham problem
is studied in \cite{Schygulla:Willmore}
for the genus 0 case and in \cite{KMR:Willmore} for arbitrary genus. Uniqueness, however,
seems to be never addressed mathematically.


It is well-known from \cite{Blaschke,MR0320956} that
the quantity $(H^2 - K) dA$ is invariant  under
M\"obius transformations,  i.e. any transformation from the group of translations (3 dimensions), rotations (3 dimensions),
uniform scalings (1 dimension) and
\textbf{\emph{sphere inversions}} (3 dimensions). If we
denote this group by $\mbox{M{\"o}b}(3)$; we have $\dim(\mbox{M{\"o}b}(3))=3+3+1+3=10$.
Here, by sphere inversion, we mean inversion about a unit sphere centered at any point in $3$-space, i.e.
\bea \label{eq:SphereInv}
i_\ba(\bx) = t_\ba \circ i \circ t_{-\ba}, \; \mbox{ where } \; i(\bx) := \frac{\bx}{\|\bx\|^2}, \quad t_\ba(\bx) = \bx+\ba.
\eea
(Sphere inversion w.r.t. a sphere with a non-unit radius can be written as one of the form \eqref{eq:SphereInv} composed with a scaling.)

The constraint functionals, namely $A$, $V$ and $M$ are only invariant under the smaller group of rigid motions ${\rm SE}(3)$.
Due to the scale-invariance of the Willmore energy,
the solution, up to homothety, of any
of the Willmore, Canham or Helfrich problems depends only on %
the
\emph{reduced volume} and \emph{reduced total mean curvature} defined by:
\bea \label{eq:v0m0}
v_0:= V_0/[(4\pi/3) (A_0/4\pi)^{3/2}], \quad m_0 := M_0/[4\pi(A_0/4\pi)^{1/2}].
\eea
This terminology is used by a group of biophysicists who have done
a plethora of  computational and physical experiments exploring the shapes of
 phospholipid vesicles. Note that $v_0$ is essentially what
a geometer would call the \emph{isoperimetric ratio}. By the isoperimetric inequality, we have
$v_0 \in (0,1]$ and $v_0=1$ is uniquely realized by a round sphere.

From now on, we think of two surfaces as the same, or that they have the same (Euclidean) shape, when they are homothetic.
By uniqueness of solution (of any one of the Helfrich, Canham or Willmore problems) we mean there is only one solution surface
up to homothety.

\subsection{Non-uniqueness in $g\geq 2$}
Given any minimizer of a Canham or Helfrich problem, one may apply to it the three dimensional family of sphere inversions \eqref{eq:SphereInv} and expect
to have enough degrees of freedom to satisfy
the reduced volume constraint or reduced volume plus mean curvature constraints, yielding a
two- or one-parameter (respectively) family of non-homothetic solutions.
This suggests that one should not expect uniqueness in general.

\gap
This hasty dimension count is easily seen to be flawed in at least specific cases. For instance,
\begin{itemize}
\item When $g=0$, the unconstrained Willmore minimizer is the round sphere and is unique, which is clearly {\bf invariant under
the whole M\"obius group}.
\item When $g=1$, the unconstrained Willmore minimizers are exactly the stereographic images into $\bR^3$
of the Clifford torus $\{[\cos u, \sin u, \cos v, \sin v]^T/\sqrt{2}: u,v \in [0,2\pi]\}$
in $\mathbb{S}^3$. For any such Clifford torus in $\bR^3$, its
Euclidean shape is {\bf invariant under 2 out of the three degrees of freedom of the sphere inversions in \eqref{eq:SphereInv}}.
(In Section~\ref{sec:StepI}, we shall establish a precise version of this fact.)
\end{itemize}
So in the first case, if we choose $v_0$ and $m_0$ to be the reduced volume and total mean curvature of the round
sphere, then the corresponding Canham or Helfrich problem must also have the round sphere as the unique solution.
In the second case, if we choose $v_0$ and $m_0$ to be the reduced volume and total mean curvature of any Clifford torus,
then we expect the corresponding Canham problem, and hence also the Helfrich problem, to have a unique solution.
The latter observation will be the focus of this article.

The dimension count, however, sounds more
convincing when the genus $g$ is 2 or above. By Hurwitz's automorphisms theorem, there can only
be a finite number -- no more than $84(g-1)$ -- conformal mappings leaving any compact genus $g$ surface invariant under homothety.
Since sphere inversions are conformal mappings,
the three-dimensional family of sphere inversions \eqref{eq:SphereInv}, when applied to any fixed compact surface of genus $g\geq 2$,
must generate a 3-dimensional family of non-homothetic surfaces.

However, uniqueness may still hold when $g\geq 2$.
To understand it better, let us first observe that instead of the 3-dimensional family of sphere inversions \eqref{eq:SphereInv},
we can instead use the 3-dimensional family of special conformal transformations
\bea \label{eq:SCT}
{\rm SCT}_\ba = i \circ t_\ba \circ i, \quad \ba \in \bR^3.
\eea
This is because for every sphere inversion $i_{\ba}$, there is a (orientation-reversing) homothety $H$
such that $i_{\ba} = H \circ {\rm SCT}_{i(\ba)}$. Moreover, since every transformation in $\mbox{M{\"o}b}(3)$
is either a homothety, an inversion, or an homothety composed with an inversion\footnote{This is a
consequence of the proof
of Liouville's theorem on conformal mappings; see, for example, \cite[Page 92]{Blair:Conformal}.}, the
non-homothetic copies of any surface $S$ under $\mbox{M{\"o}b}(3)$
can be found in
$\{ {\rm SCT}_\ba (S): \ba \in \bR^3 \} $.

So part of the (non-)uniqueness analysis boils down to the understanding of the map
$$
\bR^3 \ni \ba \stackrel{\Gamma_S}{\goto} \begin{bmatrix}
                      v({\rm SCT}_\ba (S)) \\
                      m({\rm SCT}_\ba (S)) \\
                    \end{bmatrix} \in \bR^2.
$$
Here $v()$ and $m()$ are the reduced volume and reduced total mean curvature of the argument surface; and we
 call the map $\Gamma_S$.
Being a nonlinear map, the mere fact that the co-domain has a
lower dimension than the domain does \emph{not} guarantee that
the pre-image of a given point $[v_0, m_0]^T \in {\rm Image}(\Gamma_S)$ is non-unique.
(E.g., for the map $(x_1,x_2,x_3) \mapsto x_1^2+x_2^2+x_3^2$, the pre-image of $0$ is a singleton.)
The implicit function theorem guarantees that if the differential of $\Gamma_S$ at the origin is
full rank, then indeed there is a curve through the origin, call it $\ba(t)$, such
that $\Gamma_S( \ba(t) ) = \Gamma_S ( 0 )$. To conclude, if $S$ is a particular solution
of a genus $g\geq 2$ Helfrich problem, and if ${\rm rank} (d \Gamma_S|_0) = 2$, then there must be a
one-parameter of non-homothetic solutions.

The use of special conformal transformation gives a nice expression for $d \Gamma_S|_0$:
$$
\nabla v|_0 = 6 v(0) (\mathbf{R}^A - \mathbf{R}^V), \;\;\; \nabla m|_0 = 2 m(0) (\mathbf{R}^A - \mathbf{R}^M),
$$
where $\mathbf{R}^A$, $\mathbf{R}^V$ and $\mathbf{R}^M$ are the area, volume and mean curvature centers of $S$;
see \cite[Section 5.3.1]{Seifert:Config}.
Therefore ${\rm rank} (d \Gamma_S|_0) = 2$ exactly when the three centers are not collinear. Note that the latter condition
says that $S$ must have a certain degree of asymmetry. For instance, it rules out the case when $S$
possesses 2 planes of mirror symmetry.

It is conjectured that the stereographic images of Lawson's
minimal surface $\xi_{g,1}$ in $\mathbb{S}^3$ \cite{Lawson:Minimal}
are the only $W$-minimizer of genus $g$ in $\bR^3$. The stereographic images of $\xi_{2,1}$ attain many different values
of reduced volume $v_0$ and reduced total mean curvature $m_0$. For many such values, it is observed
in \cite{LIPOWSKY:Genus2} that
the corresponding centers are not collinear and hence there is a one-parameter family of solution surfaces. (However, it is not clear if a
rigorous proof is available for this claim.)
This
non-uniqueness is called ``conformal diffusion" in the biophysics literature and is observed experimentally in a laboratory setting
\cite{MichaletBensimon:Genus2}.

\subsection{Empirical Uniqueness in genus $g=0$ and $1$}
For the genus 0 Canham problem, of which existence is shown for all $v_0\in (0,1]$ \cite{Schygulla:Willmore},
it is observed from a lot of computations (e.g. \cite{Seifert:Config,ChenYu:Helfrich}) that the solution is unique and is a surface of revolution.
When the reduced volume $v_0$ is greater than a certain value approximately equal to $0.591$, the solution surface appears to be
have an additional plane of mirror symmetry orthogonal to the axis of revolution; in this case we expect
$\mathbf{R}^A =\mathbf{R}^V$.
When $v_0$ is smaller than $0.591$, a phase transition occurs; the solution surface is a so-called  stomatocyte, which still appears to
be a surface of revolution but loses the additional plane of mirror symmetry.
When $v_0=1$, the solution is a round sphere, when $v_0\goto 0$, the solution approaches a `double sphere'.

For the genus 1 Canham problem, the existence is only established for $v_0$ in an (unknown) open interval
containing
\bea \label{eq:v0_interval}
 \big[(3/2) (2\pi^2)^{-1/4}, 1 \big);
\eea
see \cite{KMR:Willmore}. This interval is also the set of
reduced volume values attained by the M\"obius transformations of the
Clifford torus -- see Figure~\ref{fig:CliffordTori} and the next section.
The value $v_0 = (3/2) (2\pi^2)^{-1/4}$ is the reduced volume of the
surface of revolution Clifford torus
$$T_{\sqrt{2}}= \Big\{ \big[ \big(\sqrt{2} + \cos(v) \big)\cos(u), \;
\big( \sqrt{2} + \cos(v) \big)\sin(u), \; \sin(v) \big]:
u, v \in [0,2\pi] \Big \}.$$
The uniqueness of the genus 1 Canham problem on the interval \eqref{eq:v0_interval} is the focus of this paper.
When $v_0 \in (0, (3/2) (2\pi^2)^{-1/4}]$, many computations suggest that the solution surface is
unique and, similar to the genus 0 case,
is a surface of revolution; see \cite{Seifert:Config,ChenYu:Helfrich} and the references therein.

 \gap
 We therefore have the following grand conjecture:
 \begin{conjecture} \label{conjecture:Main}
The genus $g=0$ or $1$ Canham problem with any isoperimetric ratio constraint $v_0 \in (0,1)$
has a unique solution up to homothety. Moreover,
\begin{enumerate}
\item[(i)] when $g=0$, for each $v_0 \in (0,1]$ the unique solution is a surface of revolution;
\item[(ii)] when $g=1$, for each $v_0 \in \left(0, (3/2) (2\pi^2)^{-1/4} \right]$, the unique solution is a surface of revolution;
\item[(iii)] when $g=1$, for each $v_0 \in \left[(3/2) (2\pi^2)^{-1/4},1 \right)$, the unique
solution is a stereographic image into $\bR^3$ of the Clifford torus $\{[\cos u, \sin u, \cos v, \sin v]^T/\sqrt{2}: u,v \in [0,2\pi]\}$
in $\mathbb{S}^3$
or, equivalently, a
M\"obius transformation of $T_{\sqrt{2}}$.
\end{enumerate}
(When $v_0=1$ and $g=0$, it is clear that the solution is unique and is the round sphere. When $v_0=1$ and any $g\geq 1$, solution does not exist by the isoperimetric inequality.)
\end{conjecture}

An obvious difficulty in proving the uniqueness conjecture in case (i) and (ii), or uniqueness/non-uniqueness in the higher genus cases, is
that in general we do not have much information about the solutions of the Canham or Helfrich problems.
As a starting point, we explore the third case of Conjecture~\ref{conjecture:Main}, which appears to be the most tractable.

\subsection{This paper}
%
\begin{figure}
\centerline{
\begin{tabular}{ccccc}
 \hspace{10pt}{\small $a=0$} & {\small $\cdots$} & {\small $\rightarrow$} & {\small $\cdots$} & \hspace{4pt}{\small $a\approx\sqrt{2}-1$}\vspace{-3pt} \\
 \hspace{10pt}{\small $v_0=0.71$} & {\small $v_0=0.75$} & {\small $v_0=0.85$} & {\small $v_0=0.95$} & \hspace{4pt}{\small $v_0=0.99$} \vspace{-1pt} \\
 \includegraphics[angle=90,origin=c,height=1.2in]{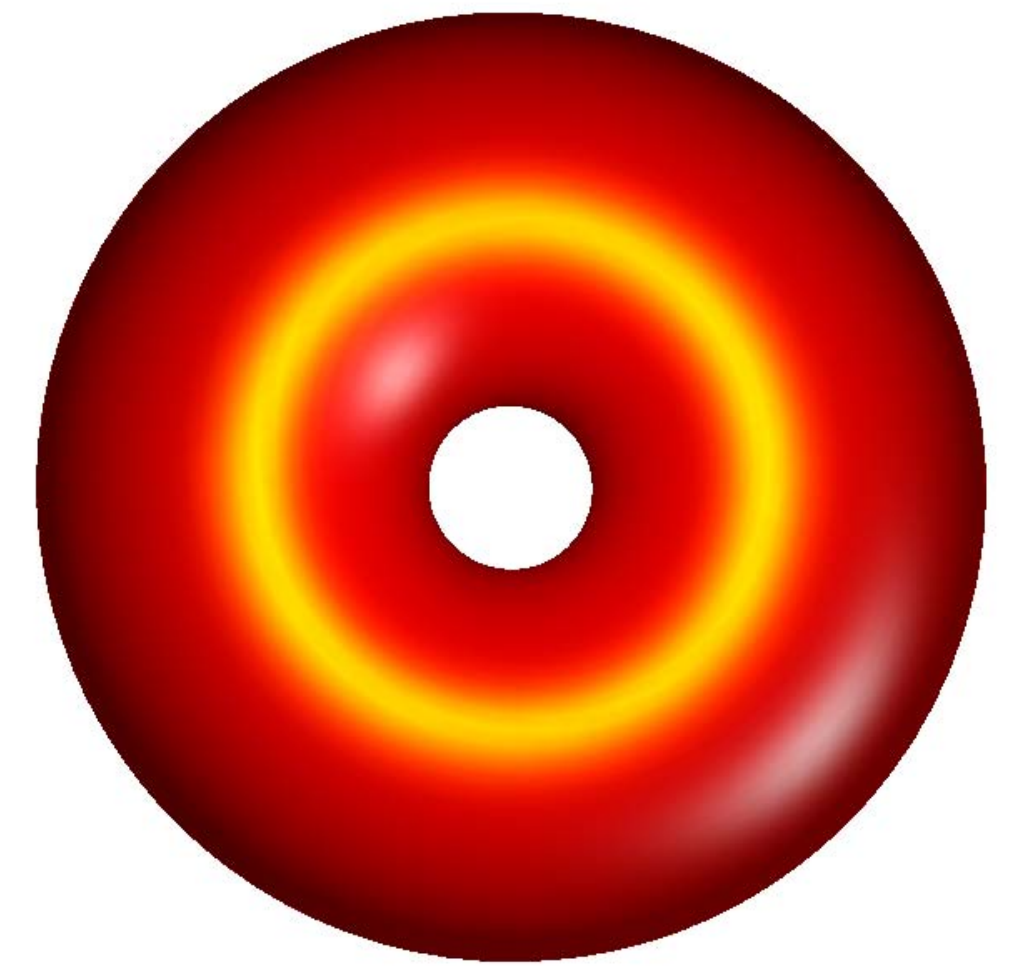} \hspace{-.1in}
 &\includegraphics[angle=90,origin=c,height=1.2in]{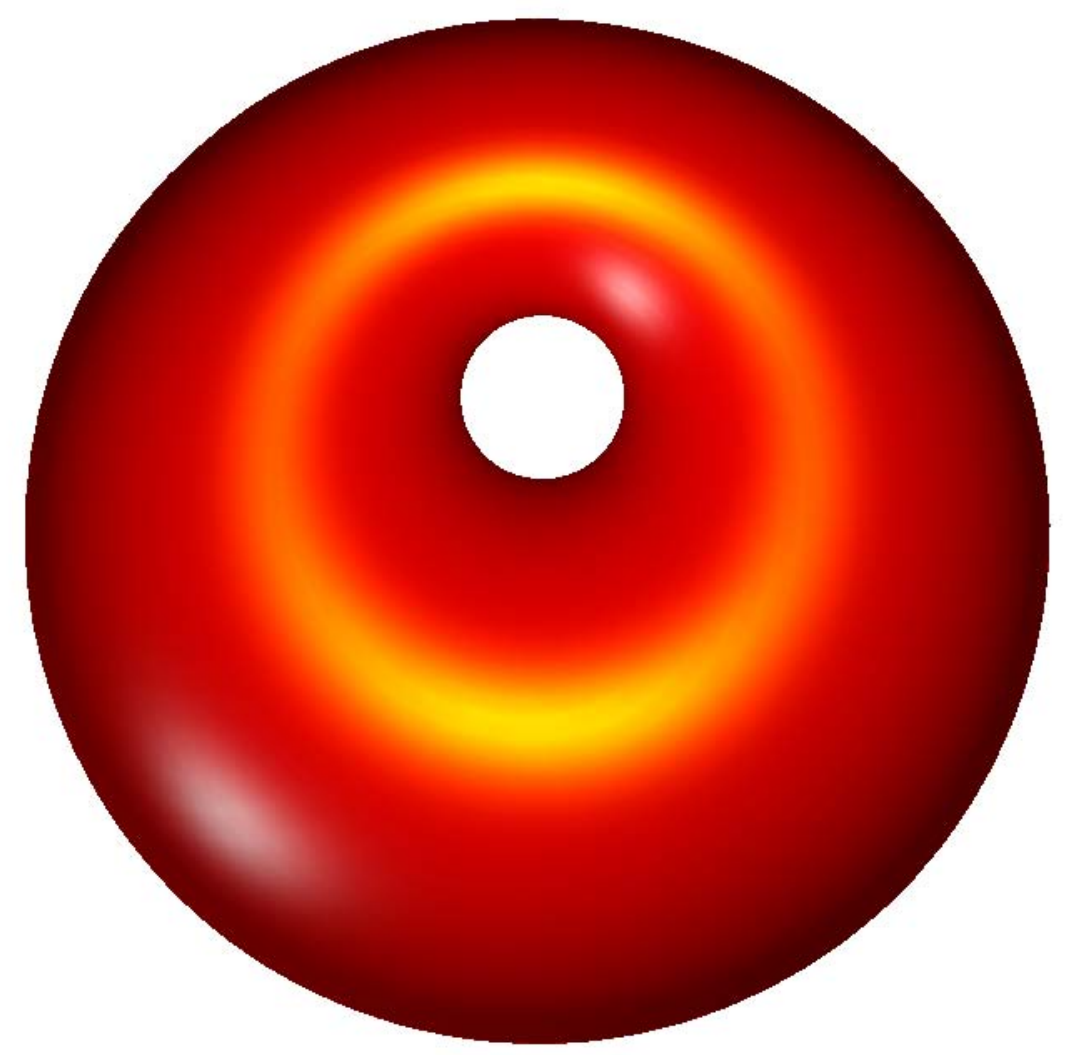}\hspace{-.1in}
 &\includegraphics[angle=90,origin=c,height=1.2in]{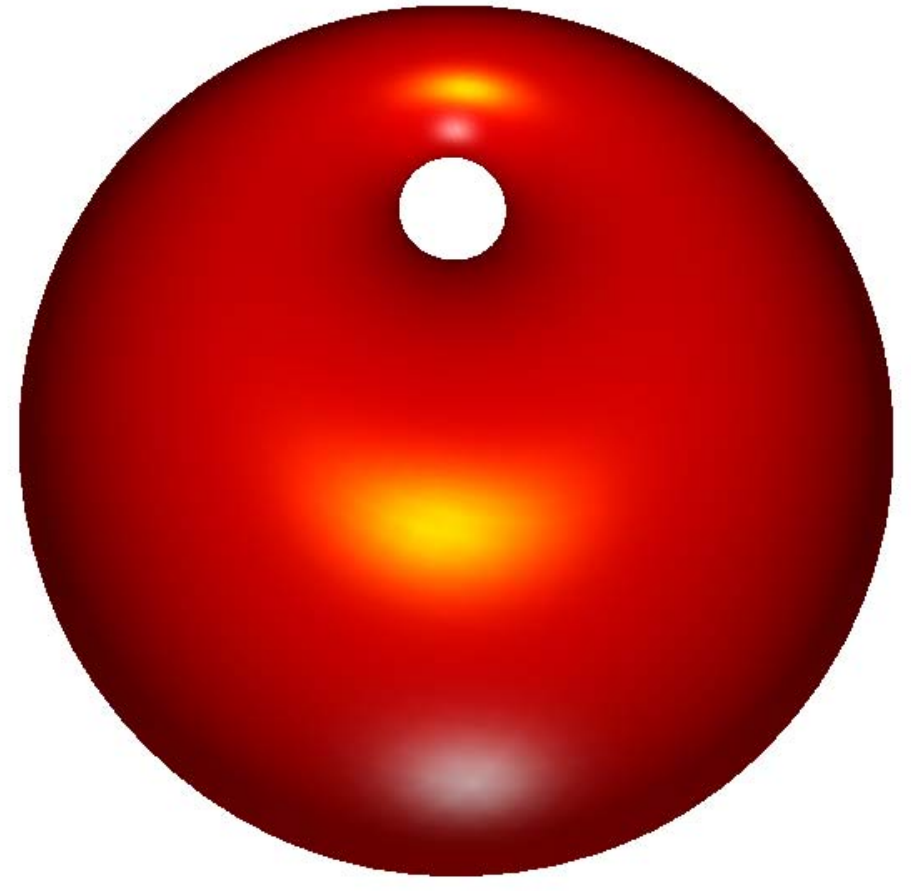}\hspace{-.1in}
 &\includegraphics[angle=90,origin=c,height=1.2in]{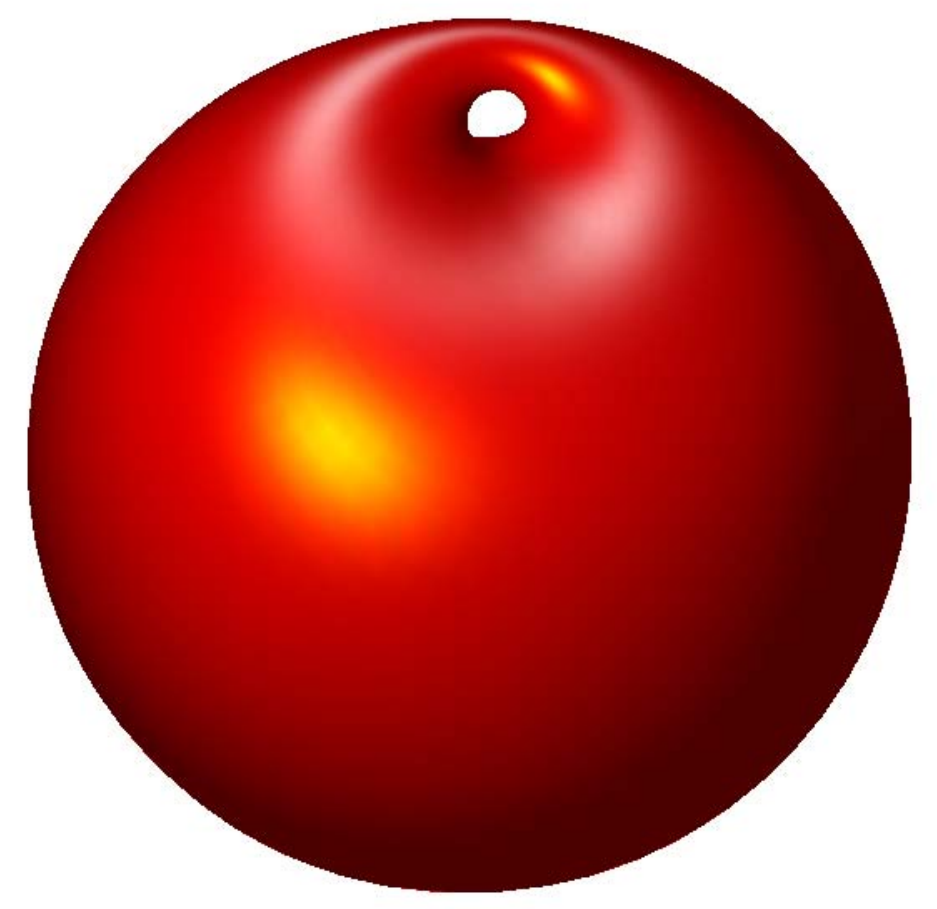}\hspace{-.1in}
 &\includegraphics[angle=90,origin=c,height=1.2in]{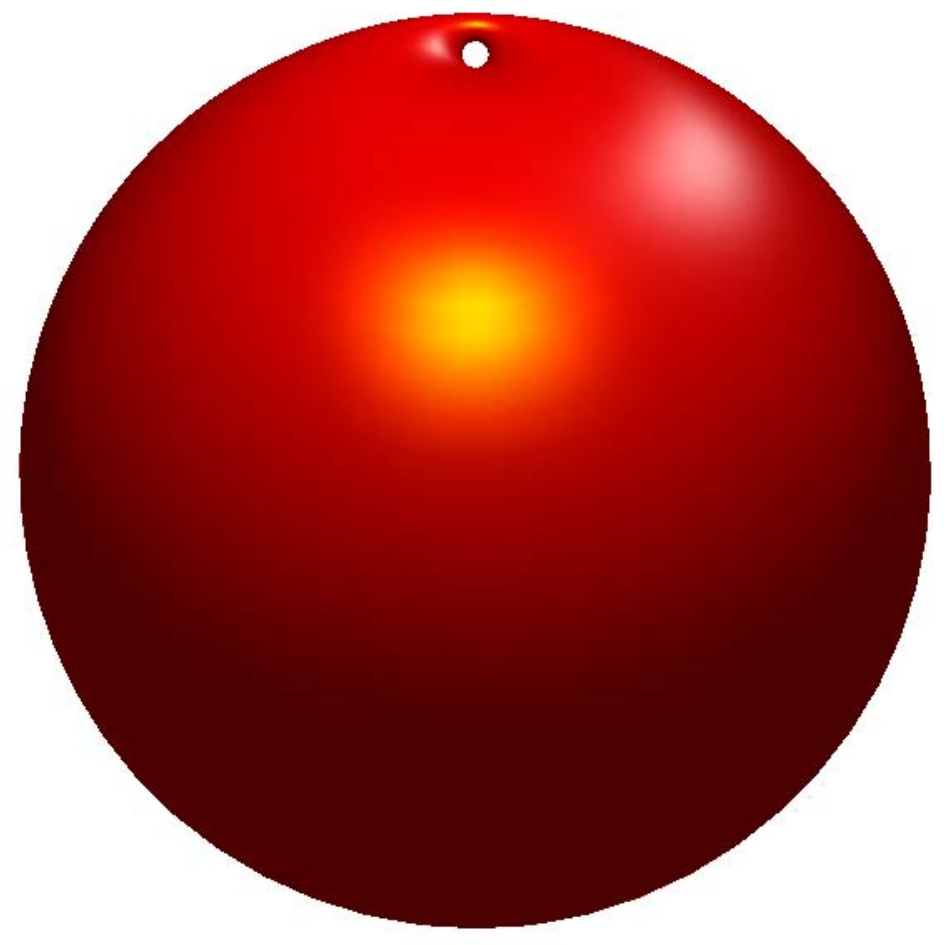}\hspace{-.1in}
\end{tabular}}
\caption{$i_{(a,0,0)}(T_{\sqrt{2}})$ as $a$ increases from $0$ to $\sqrt{2}-1$. By Theorem~\ref{thm:MainStepI},
these are all the possible non-homothetic
images of $T_{\sqrt{2}}$ under $\mbox{M\"ob}(3)$.
By Pappus's centroid theorems, the reduced volume of the surface of revolution Clifford torus
$T_{\sqrt{2}}$ is $(3/2) (2\pi^2)^{-1/4} \approx 0.71$.
Notice that $i_{(0,0,0)}(T_{\sqrt{2}})=T_{\sqrt{2}}$, viewing $T_{\sqrt{2}}$ as a point set.
(As an oriented surface, $T_{\sqrt{2}}$ is turned inside out by $i=i_{(0,0,0)}$.)
When $a$ approaches $\sqrt{2}-1$,
$i_{(a,0,0)}(T_{\sqrt{2}})$
approaches a round sphere, and hence
should have a reduced volume close to $1$; see Section~\ref{sec:Step2} for a proof.}\label{fig:CliffordTori}
\end{figure}

\gap
To establish Conjecture~\ref{conjecture:Main}(iii), we propose the following four steps:
\begin{itemize}
\item[I.] Prove that the set of all non-homethetic images of $T_{\sqrt{2}}$ under $\mbox{M\"ob}(3)$ corresponds
 exactly to the one-parameter family
 \bea \label{eq:cyclides}
 \{i_{(a,0,0)}(T_{\sqrt{2}}):  a\in [0,\sqrt{2}-1)\}.
 \eea
 In other words, the cyclides
depicted in Figure~\ref{fig:CliffordTori} are exactly the set of all non-homothetic Clifford tori.
 This is established in Theorem~\ref{thm:MainStepI} of Section~\ref{sec:StepI}.

\item[II.] With this result,
 the conjecture follows if we can show:
\bea
\Iso: [0,\sqrt{2}-1) \goto [(3/2) (2\pi^2)^{-1/4},1), \quad \Iso(a) := v\big( i_{(a,0,0)}(T_{\sqrt{2}}) \big)
\eea
 is a bijection.
If so, then each $v_0 \in [(3/2) (2\pi^2)^{-1/4},1)$ corresponds
 to one and only one Clifford torus, namely $i_{(\Iso^{-1}(v_0),0,0)}(T_{\sqrt{2}})$, with isoperimetric ratio $v_0$, which
 must be the unique solution
 of the genus 1 Canham problem with $v_0$ as the constrained isoperimetric ratio.

To prove that $\Iso$ is a bijection, it suffices to show that
$\Iso$ is monotonic increasing and
\bea \label{eq:lim1}
\lim_{a \goto \sqrt{2}-1} \Iso(a) = 1.
\eea
In Section~\ref{sec:Step2}, we establish Theorem~\ref{thm:Rounding}, which is a more general version of \eqref{eq:lim1}.

\item[III.]  To prove that $\Iso$ is monotonic increasing, we venture into the realm of
special functions. We make the observation that the area and enclosing volume of the cyclides in \eqref{eq:cyclides},
denoted by $A(a)$ and $V(a)$, can be extended analytically to the disc $\{z: |z|<\sqrt{2}-1\}$ on the complex plane. Moreover,
the coefficients $(a_n)_{n\geq0}$ and $(v_n)_{n\geq 0}$  of their power series at $z=0$ are holonomic, or P-recursive, sequences,
i.e. they satisfy linear recurrences with polynomial coefficients. We work out explicitly these P-recurrences
in Section~\ref{sec:StepII}.

Since $\Iso^2(a)/(36\pi)  = V^2(a)/A^3(a)$, $\Iso$ is monotonic increasing iff the logarithm of the
right-hand side is. But then we have
$$
\frac{d}{da} \ln \frac{V^2(a)}{A^3(a)}  = \frac{2 V'(a) A(a) - 3 V(a) A'(a)}{V(a) A(a)},
$$
so $\Iso$ is monotonic increasing iff $2 V'(a) A(a) - 3 V(a) A'(a)>0$ on $[0,\sqrt{2}-1)$.
The fact that $A(z)$ and $V(z)$ are holonomic implies
that $D(z):= 2 V'(z) A(z) - 3 V(z) A'(z)$ is also holonomic; the coefficients $(d_n)_{n\geq0}$
of the power series of $D(z)$ at $z=0$ follows the P-recurrence \eqref{eq:d_recur} derived in Section~\ref{sec:StepII}.

The monotonicity of $\Iso$ follows if all
the terms defined by the P-recurrence \eqref{eq:d_recur} are positive.
\item[*IV.]
Prove that all terms defined by the P-recurrence \eqref{eq:d_recur} are positive.

This last step is out of the scope of this paper. It is well-known to experts in holonomic functions that positivity of a P-recurrence is
difficult to establish when its characteristics polynomial
has a non-simple dominant root, as is the case of
\eqref{eq:d_recur}. We can, however, use the existing tools to check that
the sequence is \emph{eventually positive}; see Section~\ref{sec:StepII} and the remarks in Section~\ref{sec:Final}.
\end{itemize}

To summarize, the result of this article is:
\begin{proposition}
Assuming the positivity of the P-recurrence \eqref{eq:d_recur},
Conjecture~\ref{conjecture:Main}(iii) holds, i.e. for every isoperimetric ratio
$v_0 \in \left[(3/2) (2\pi^2)^{-1/4},1 \right)$,
there is a unique M\"obius transformation of $T_{\sqrt{2}}$ with isoperimetric ratio $v_0$.
\end{proposition}

Steps I-III are carried out in the next three sections.

\section{Step I: Non-homothetic Clifford tori} \label{sec:StepI}

Let $T_R := \left\{
(R+ \cos v) \cos u,(R+ \cos v) \sin u, \sin v): u,v \in [0,2\pi] \right\}$, a torus
with major radius $R \in (1,\infty)$, minor radius $1$, and the $z$-axis as the axis of revolution.
Let $i_{(x,y,z)}$ be the
inversion map about the unit sphere centered at $(x,y,z)$ of $\bR^3$. Our
goal is to characterize all the Euclidean shapes of the Clifford tori, i.e.
 we would like to find a parametrization of the `shape space'
\bea \label{eq:ShapeSpace_Clifford}
\left\{i_{(x,y, z)}(T_{\sqrt{2}}): (x,y,z) \in \bR^3 \backslash T_{\sqrt{2}} \right\} \Big/ {\rm Hom}(3).
\eea
Here `$/ {\rm Hom}(3)$' means we identify two point sets if they can be transformed from one to another by a
homothety in $\bR^3$. Since we are primarily interested in Euclidean shapes here, we avoid sphere inversions centered
at points on $T_R$ itself. To help us gain a better understanding of the underlying structure, we also study the more general
shape space
\bea \label{eq:ShapeSpace}
\left\{i_{(x,y, z)}(T_R): R>1, \; (x,y,z) \in \bR^3 \backslash T_R \right\} \Big/ {\rm Hom}(3).
\eea

\gap
\noindent {\bf Maxwell's characterization of a cyclide.}
It is well-known that any (torodial) cyclide $\mathfrak{C}$ has two orthogonal planes of mirror symmetry;
see, for example, \cite{Maxwell:Cyclide,Boehm:cyclide,CDH:cyclide}.
We make the observation that the Euclidean shape of a toroidal cyclide $\mathfrak{C}$
is uniquely determined by certain measurements of
the cross section of $\mathfrak{C}$ with either one of the two symmetry planes.

\begin{figure}[ht]
\centering
\begin{tabular}{cccc}
\hspace{-.2in}
\includegraphics[height=1.3in]{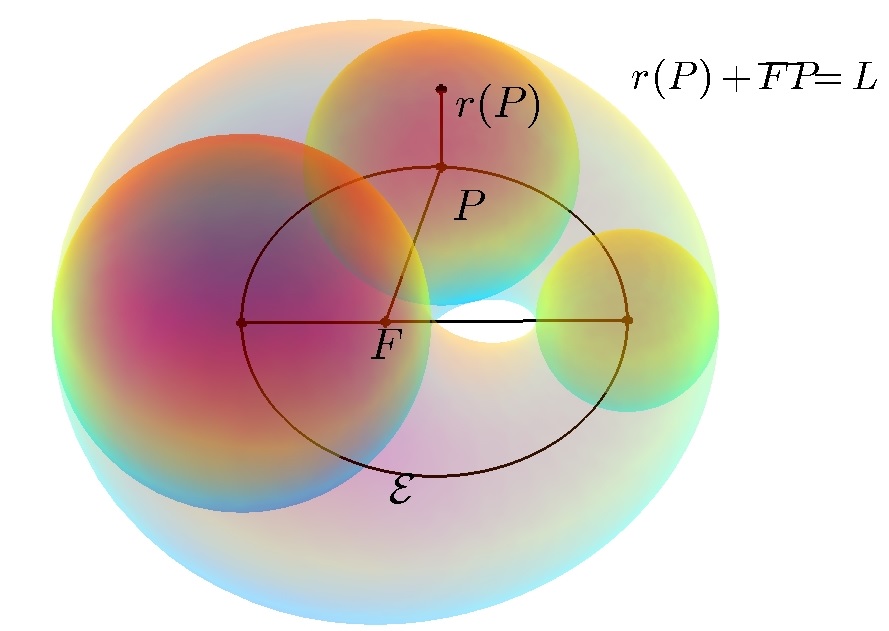}
& \hspace{-.2in}
\includegraphics[height=1.3in]{CyclideP1P2}
& \hspace{-.2in}
\includegraphics[height=1.1in]{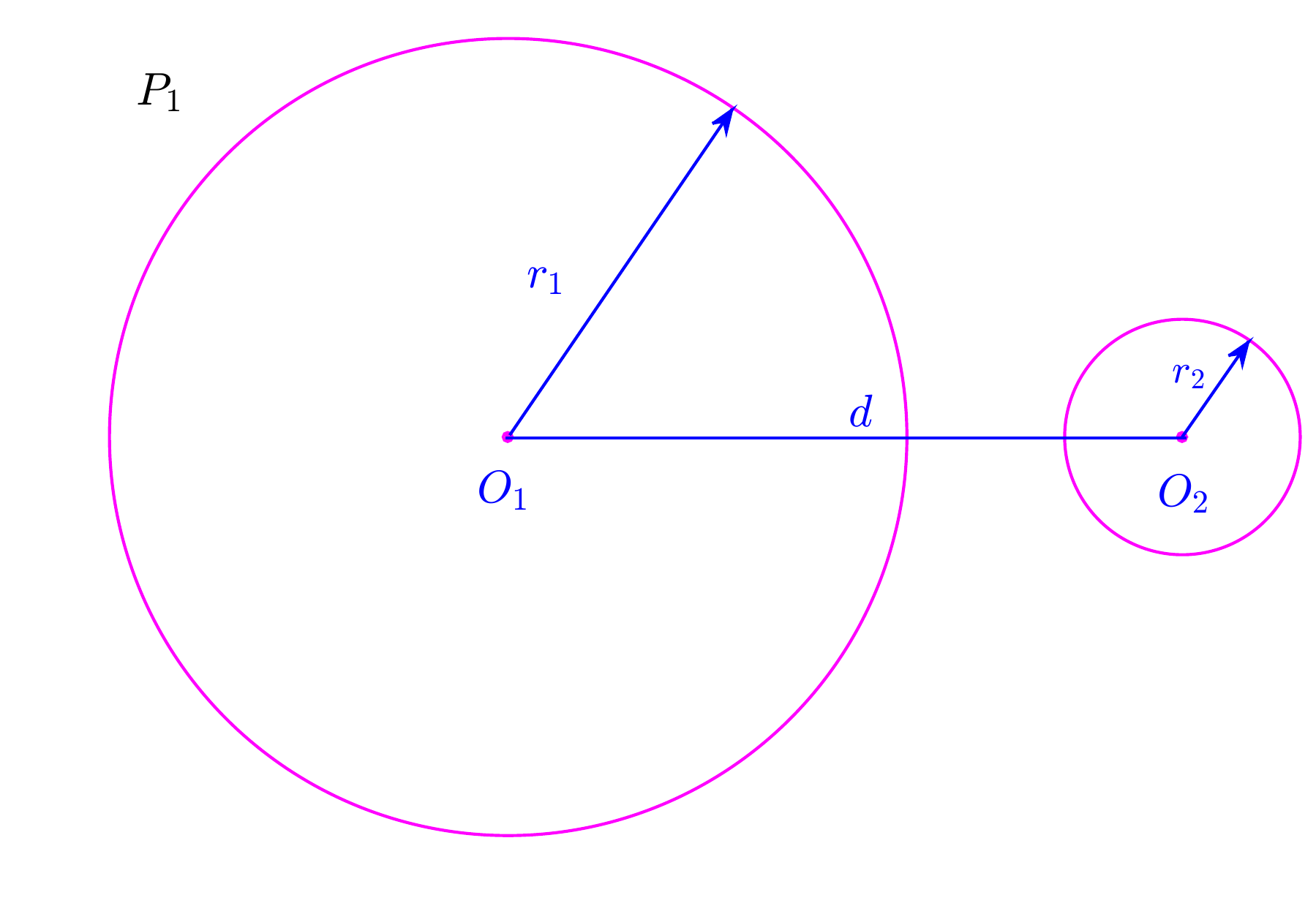}
& \hspace{-.2in}
\includegraphics[height=1.1in]{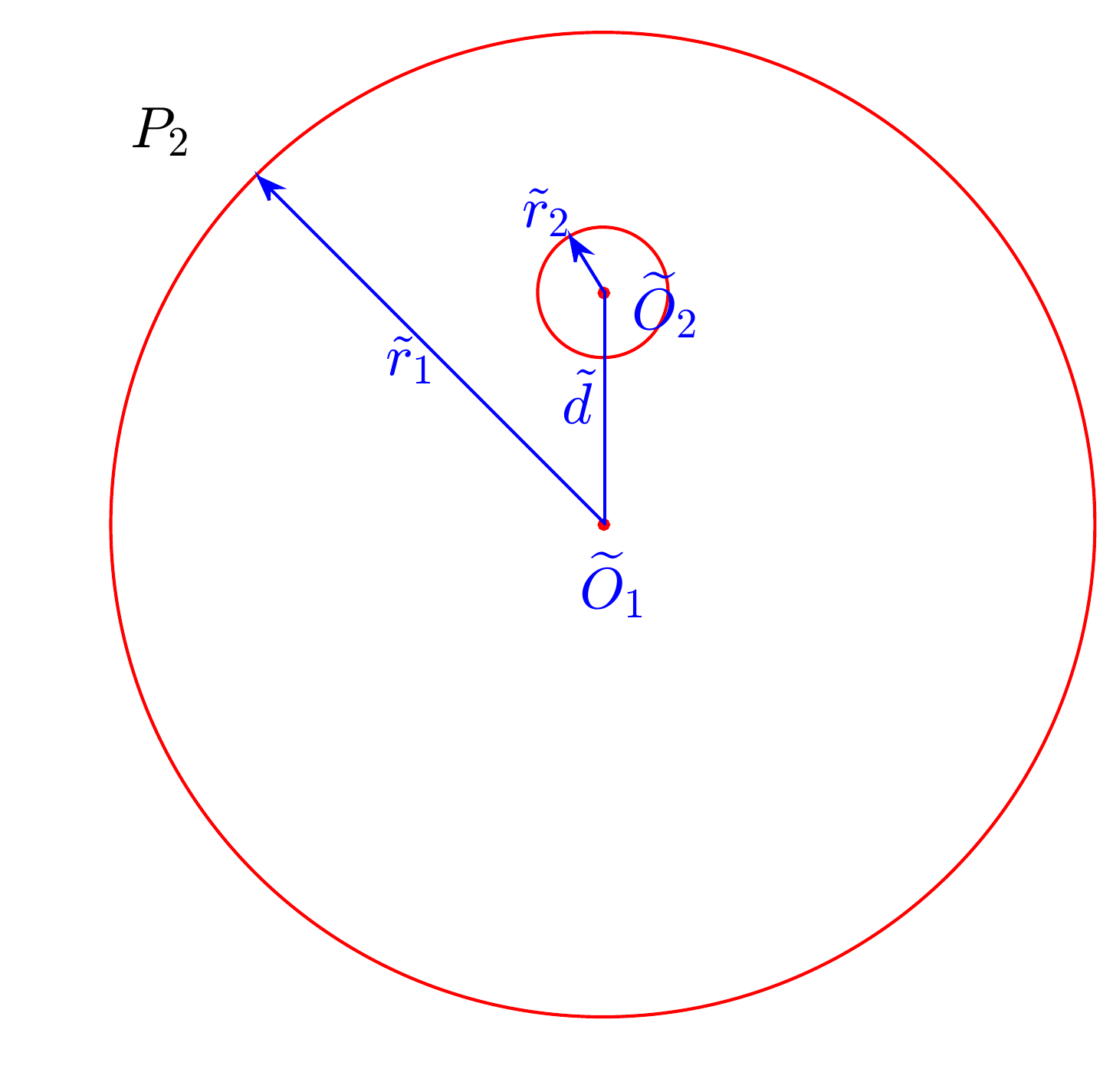}\\
\hspace{-.2in}
(a)
& \hspace{-.2in}
(b)
& \hspace{-.2in}
(c) $i_{(x,y, z)}(T_R) \cap P_1$
& \hspace{-.2in}
(d) $i_{(x,y, z)}(T_R) \cap P_2 $
\end{tabular}
\caption{(a) Maxwell's characterization of a torodial cyclide, (b) Two planes of mirror symmetry,
(c)-(d) Cross sections of $i_{(x,y, z)}(T_R)$ with $P_1$ and $P_2$
}
\label{fig:MaxwellCyclide}
\label{fig:CyclideP1P2}
\label{fig:CrossSections}
\end{figure}

We use Maxwell's characterization of cyclides \cite{Maxwell:Cyclide,Boehm:cyclide,CDH:cyclide}:
any cyclide $\mathfrak{C}$
is the envelope of all the spheres centered at the points $P$ on a given ellipse $\mathcal{E}$
with radii $r(P)$, $P \in \mathcal{E}$, satisfying
$r(P) + \overline{FP} = L$, where $F$ is one of the foci of $\mathcal{E}$ and $L$ is a constant in a suitable range.
We can think of $L$ as the length of a taut
string attached in one end to $F$; the string slides smoothly on $\mathcal{E}$ and
traces out spheres with the other end. See Figure~\ref{fig:MaxwellCyclide}. Under this characterization,
$\mathfrak{C}$ is a torodial cyclide if and only if
$$
a > L-a > f,
$$
where $a$, $f$ and $L$ are the major radius of $\mathcal{E}$, the focal length of $\mathcal{E}$,
and the length of the string, respectively. Moreover,
the Euclidean shape of $\mathfrak{C}$ can be characterized by the ratio
$a:f:L$.\footnote{This already explains why the shape space \eqref{eq:ShapeSpace}
is two-dimensional.}

The major axis of $\mathcal{E}$ lies on the intersecting line of the two symmetry planes of $\mathfrak{C}$.
In the following, $P_1$ refers to the symmetry plane where $\mathcal{E}$ lies,
whereas $P_2$ ($\bot P_1$) refers to the other symmetry plane.
The cross section $\mathfrak{C} \cap P_1$ consists of two circles {\bf exterior to each other},
whereas the cross section $\mathfrak{C} \cap P_2$ consists of two circles
with {\bf one lying inside the other} (see Figure~\ref{fig:CrossSections}).

Denote the radii of the two circles in $\mathfrak{C} \cap P_1$ by $r_1$ and $r_2$ and the distance between the two centers by
$d$ (see Figure~\ref{fig:CrossSections}). Similarly, let $\tilde{r}_1$ and $\tilde{r}_2$ be the radii of the two
circles in $\mathfrak{C} \cap P_2$ and $\tilde{d}$ be the distance between the two centers.
By convention, $r_1 \geq r_2$, $\tilde{r}_1 \geq \tilde{r}_2$.  The three sets of measurements $(r_1, r_2, d)$, $(\tilde{r}_1, \tilde{r}_2, \tilde{d})$ and $(a,f,L)$ of a cyclide $\mathfrak{C}$ are related by the following equations:
\beqn 
a = \frac{d}{2}, \;\; f = \frac{r_1-r_2}{2}, \;\; L = \frac{d+r_1+r_2}{2}.
\eeqn
\beqn
\label{eq:P1toP2}
\tilde{r}_1 = \frac{d + (r_1+r_2)}{2}, \;\;
\tilde{r}_2 = \frac{d - (r_1+r_2)}{2}, \;\;
\tilde{d} = r_1 - r_2.
\eeqn
Since the maps $(a,f,L) \mapsto (r_1,r_2,d)$ and $(r_1,r_2,d) \mapsto (\tilde{r}_1, \tilde{r}_2, \tilde{d})$
are linear isomorphisms,
we conclude that:
\begin{lemma} \label{lemma:Maxwell}
Each of the three ratios
$$
a:f:L, \quad r_1:r_2:d \;\;\mbox{and } \;\; \tilde{r}_1:\tilde{r}_2:\tilde{d}
$$
determines
the Euclidean shape of the cyclide $\mathfrak{C}$.
\end{lemma}

\begin{figure}[ht]
\centerline{
\begin{tabular}{c}
\includegraphics[height=2.8in]{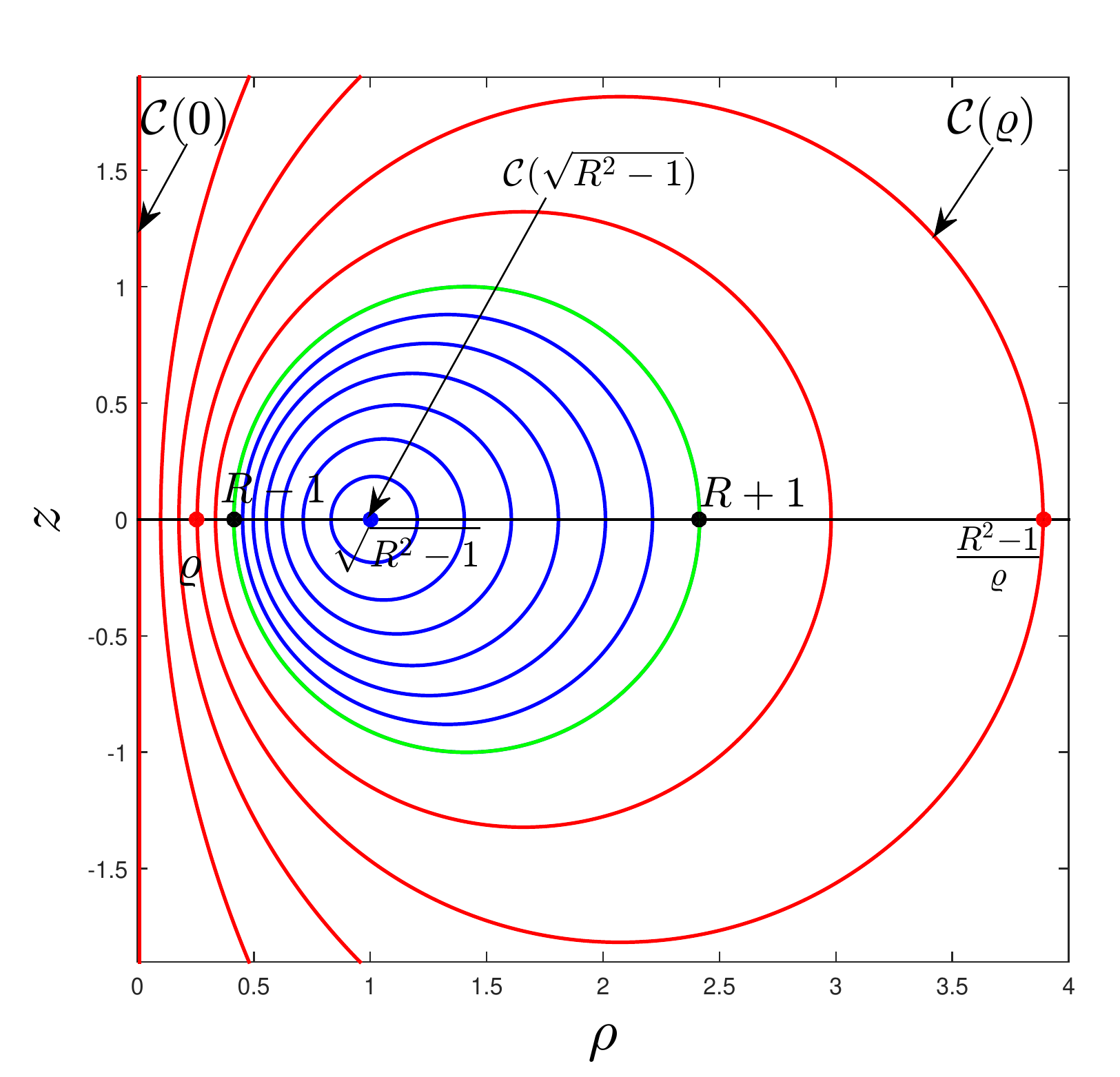}
\end{tabular}}
\caption{
$\mathcal{C}(\varrho)$ for various values of $\varrho \in (0,R-1)$ (in red) and
$\varrho \in (R-1,\sqrt{R^2-1})$ (in blue). Note that $\mathcal{C}(\varrho) = \mathcal{C}((R^2-1)/\varrho)$
and $\mathcal{C}(\sqrt{R^2-1})$ degenerates into a point.
}
\label{fig:Circles}
\end{figure}

For any $\varrho > 0$, let $\mathcal{C}(\varrho) = \mathcal{C}(\varrho; R)$ be the circle in the $\varrho$-$z$ plane with
a diameter connecting $(\varrho,0)$ and $\left( (R^2-1)/\varrho,0 \right)$; see Figure~\ref{fig:Circles}.
By convention, $\mathcal{C}(0) = \mathcal{C}(\infty) $ is the $z$-axis.
In general, we have
$$
C(\varrho) = C( (R^2-1)/\varrho ).
$$

These circles on the plane can be extended to the following tori in 3-D:
\bea \label{eq:Trho}
\mathcal{T}(\varrho) := \mathcal{T}(\varrho; R) := \left\{ (\rho \cos(\theta), \rho \sin(\theta), z): (\rho,z) \in \mathcal{C}(\varrho), \; \theta \in [0,2\pi] \right\}.
\eea
For any fixed $R$, the torus
$\mathcal{T}(\varrho)$ lies
completely outside, on, or inside the torus $T$
when $\varrho \in [0,R-1) \cup (R+1, \infty]$, $\varrho = R\mp 1$, or $\varrho \in (R-1, R+1)$, respectively. In particular,
$\mathcal{T}(R\pm 1; R) = T_R$.
On the $\rho$-$z$ plane, these correspond to the red, green and blue circles in Figure~\ref{fig:Circles}.
While the one-parameter family of circles
$$
\left\{
\mathcal{C}(\varrho): \varrho \in [0,\sqrt{R^2-1}]
\right\}
$$
partitions the $\rho$-$z$ plane,\footnote{Any $(\rho, z)$, $\rho>0$, lies on the circle $\mathcal{C}(\varrho^+) = \mathcal{C}(\varrho^-)$, where
$$
\varrho^\pm = \frac{\left(\rho^2+z^2+R^2-1 \right) \pm
\sqrt{\left(\rho^2+z^2+R^2-1 \right)^2 - 4\rho^2(R^2-1)}}{ 2\rho}.
$$}
the corresponding one-parameter family of tori
$$
\left\{ \mathcal{T}(\varrho): \varrho \in [0,\sqrt{R^2-1}] \right\}
$$
partitions $\bR^3$. We shall see that how these circles and tori characterize the shape spaces
\eqref{eq:ShapeSpace_Clifford} and \eqref{eq:ShapeSpace}.

\begin{theorem} \label{thm:FirstHomothetic}
For any fixed $R \in (1,\infty)$ and $\varrho \in [0,\infty] \backslash\{R-1,R+1\}$, all the cyclides in
\bea
\label{eq:SameShape}
\Big\{i_{(x,y, z)}(T_R): \quad (x,y,z) \in \mathcal{T}(\varrho; R) \Big\},
\eea
are homothetic in $\bR^3$.
\end{theorem}
\pf
We divide the proof into 3 steps.

\gap
\noindent
$1^\circ$ By rotational symmetry, the shape of $i_{(\rho \cos(\theta), \rho \sin(\theta), z)}(T_R)$ is independent of $\theta$.
So it suffices to prove that all cyclides of the form
$$i_{(\rho, 0, z)}(T_R), \quad (\rho, z) \in \mathcal{C}(\varrho),$$
are homothetic.

By Lemma~\ref{lemma:Maxwell},
the Euclidean shape of $i_{(\rho,0,z)}(T_R)$ is determined by the measurements
of its cross section at the $x$-$z$ plane. Denote by $P$ the $x$-$z$ plane and $\pi: \bR^3 \goto \bR^2$ be the ortho-projection
map onto $P$. Note that
\bea \label{eq:Cylcides_P}
\pi\left( i_{(\rho,0,z)}(T_R) \cap P \right) = i_{(\rho,z)} \left( \pi(T_R \cap P) \right).
\eea
Here $i_{(\rho, z)}$ stands for the circle
inversion map in 2-D with respect to the
unit circle centered at $(\rho,z)$.
Note that $P$ is a symmetry plane of the cyclides \eqref{eq:Cylcides_P} and that the cross section
\eqref{eq:Cylcides_P} consists of a circle pair.
Therefore, by (the implication of) Lemma~\ref{lemma:Maxwell}, it suffices to check that these circle pairs corresponding to
different $(\rho,z) \in \mathcal{C}(\varrho)$ are all homothetic.
We have reduced the problem into one of plane geometry.

\begin{figure}[ht]
\centerline{
\begin{tabular}{cc}
\includegraphics[height=1.8in]{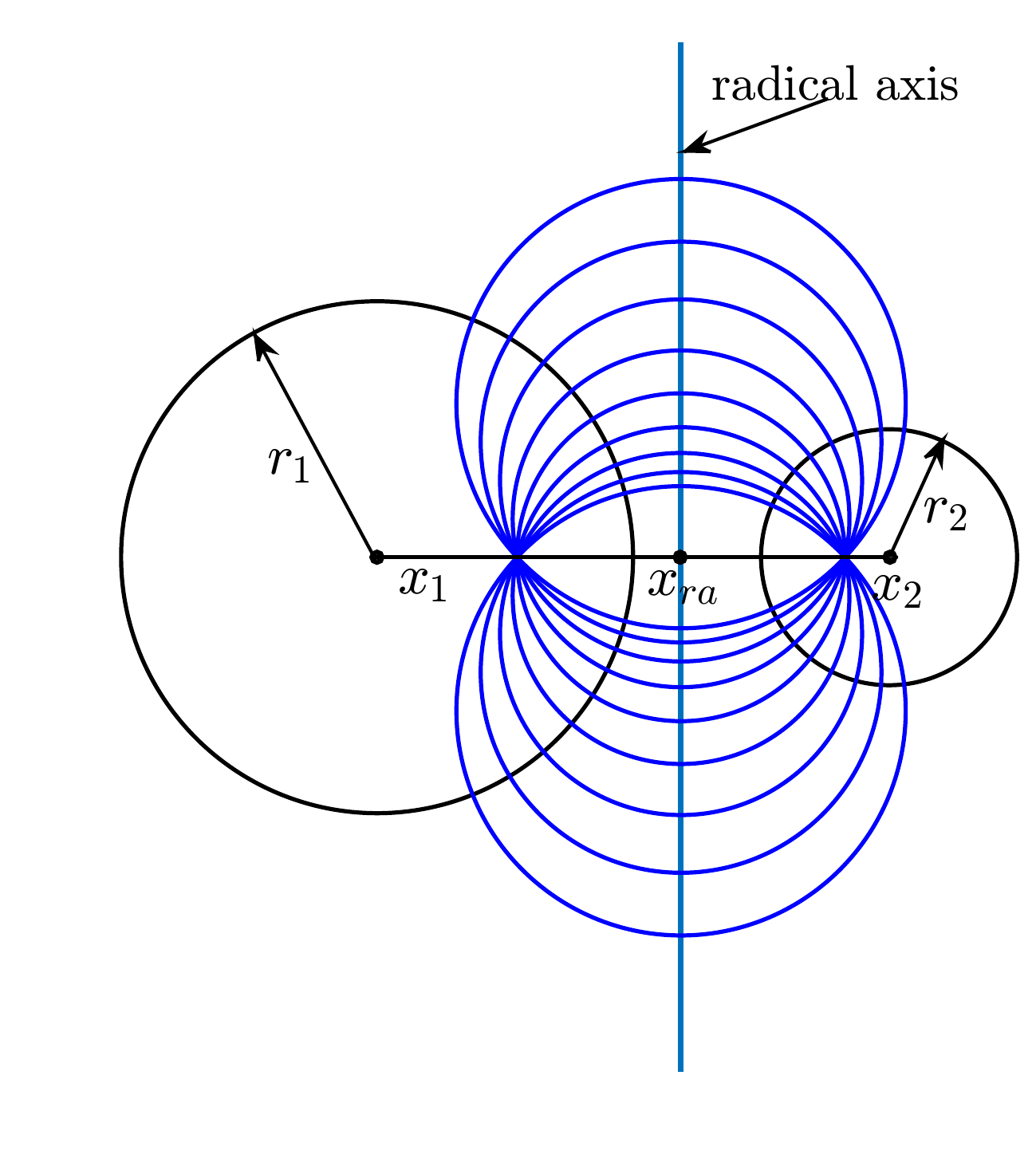} &
\includegraphics[height=1.8in]{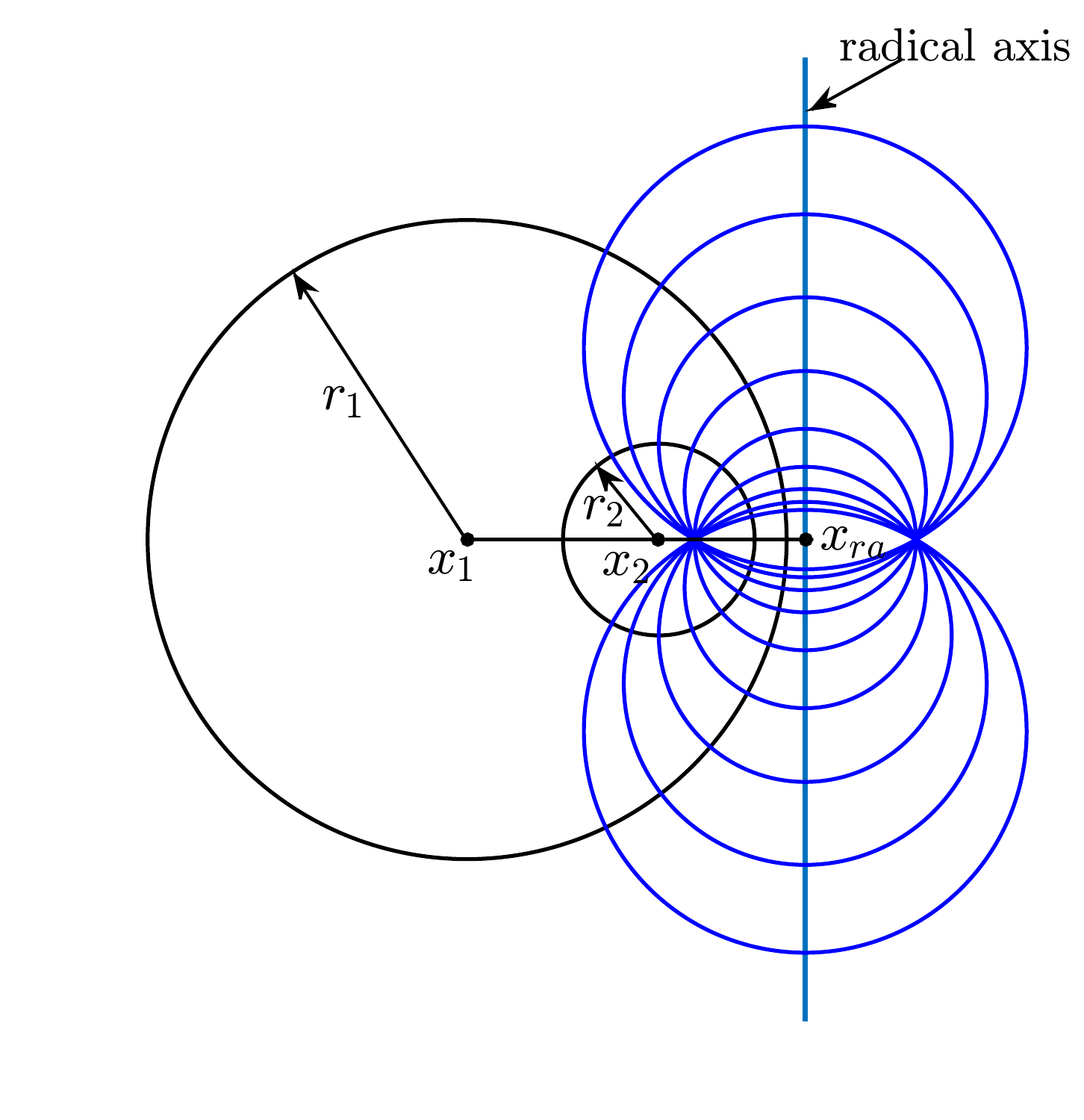} \\
\end{tabular}}
\caption{Radical Axis of a circle pair: (left) two circles exterior to each other; (right) one circle lying inside the other.
The blue circles meet the circle pair orthogonally.}
\label{fig:RadicalAxis}
\end{figure}

\gap
\noindent
$2^\circ$ We recall a well-known fact about circle inversion. If we invert two circles centered
at $(x_1, 0)$ and $(x_2,0)$ with radii $r_1$ and $r_2$ about a circle centered anywhere on the line
\bea \label{eq:RadicalAxis}
\left \{
(x_{ra},y) \; \Big| \; x_{ra} = \frac{(x_
2^2 - x_1^2) + (r_1^2-r_2^2)}{2(x_2- x_1)} \right \},
\eea
the resulting circle pair
is homothetic to the original circle pair. This line is called the {\it radical axis} of the circle pair;
see Figure~\ref{fig:RadicalAxis}.

We first determine the image of the circle pair
$\pi(T_R \cap P)$ under the circle inversion $i_{(\varrho, 0)}$.
The circle pairs in $\pi(T_R \cap P)$ consist of two unit circles with diameters $A_1 B_1$ and $A_2 B_2$, both on the $x$-axis,
with $A_1=(R-1,0)$,
$B_1=(R+1,0)$, $A_2=(-(R+1),0)$ and $B_2=(-(R-1),0)$.
The images of $A_1$, $B_1$, $A_2$, $B_2$ under $i_{(\varrho, 0)}$, denoted by $A_1'$, $B_1'$, $A_2'$, $B_2'$, again lie on the $x$-axis
and form the diameters $A_1'B_1'$, $A_2'B_2'$ of circle pair in $i_{(\varrho, 0)}( \pi(T_R \cap P) )$.
\begin{itemize}
\item When $\varrho \in (0, R-1)$, $B_2'<A_2'<B_1'<A_1'$.
\item When $\varrho \in (R+1,\infty)$, $B_1'<A_1'<B_2'<A_2'$.
\item When $\varrho \in (R-1,R+1)$, $A_1'<B_2'<A_2'<B_1'$.
\end{itemize}\footnote{Here and below, $A<B$ simply means $A$ is on the left of $B$ for
 two points $A$ and $B$ are on the first axis of $\bR^2$.}
In the first two cases, the circle pair are exterior of each other, as in Figure~\ref{fig:RadicalAxis}(a);
in the last case, one circle lies inside the other, as in Figure~\ref{fig:RadicalAxis}(b).
In any case, the resulting
circle pair has the following radii and centers:
\begin{equation} \label{eq:cyclideCrossSections}
\begin{aligned}
\mathfrak{r}_1  = \frac{|A_1'-B_1'|}{2} = \frac{1}{|(\varrho - R)^2-1|}, &
\;\;\;
\mathfrak{r}_2  = \frac{|A_2'-B_2'|}{2} = \frac{1}{(\varrho + R)^2-1} \\
O_1  = \frac{A_1'+B_1'}{2} =\Big(\varrho - \frac{\varrho - R}{(\varrho - R)^2 -1}, 0 \Big), &
\;\;\;
O_2 = \frac{A_2'+B_2'}{2} =\Big(\varrho - \frac{\varrho + R}{(\varrho + R)^2 -1}, 0 \Big).
\end{aligned}
\end{equation}
By \eqref{eq:RadicalAxis} and \eqref{eq:cyclideCrossSections}, the radical axis of the circle pair $i_{(\varrho, 0)} \left( \pi(T_R \cap P) \right)$
is given by $\{(\rho_{ra},z): z \in \bR\}$ where
$$\rho_{ra} = \varrho - \frac{\varrho}{\varrho^2+1-R^2}.$$

Now the circle pairs in
\bea \label{eq:CirclePairs}
\Big\{i_{(\rho_{ra},z)}\circ i_{(\varrho,0)} ( \pi(T_R \cap P) ): \; z \in \bR\Big\}
\eea
are all homothetic. The theorem is proved if we
show that every circle pair in
$\Big\{i_{(\rho, z)}(\pi( T_R \cap P )): \; (\rho,z) \in \mathcal{C}(\varrho) \Big\}$
is homothetic to some circle pair in \eqref{eq:CirclePairs}.
We do so in the last step of the proof.

\gap
\noindent
$3^\circ$ Since an arbitrary composition of inversions can be written as a
composition of an inversion (of radius 1) with a homothety (see \cite[Page 92]{Blair:Conformal}),
\bea \label{eq:InvComp}
i_{(\rho_{ra},z)} \circ i_{(\varrho,0)} = \mathcal{H} \circ i_{(\rho_1,z_1)}.
\eea

We can determine $(\rho_1,z_1)$ using the following properties of an inversion $i_O$ to find $(\varrho_1,z_1)$:
$i_O (O) = \infty$,
$i_O(\infty) = O$,
and $i_O(Q_1) = Q_2 \Leftrightarrow i(Q_2) = Q_1$.
By the first property,
$$i_{(\rho_{ra},z)} \circ i_{(\varrho,0)}(\rho_1,z_1) = \mathcal{H} \circ i_{(\rho_1,z_1)} (\rho_1,z_1) = \infty.$$
By the second property,
$$i_{(\varrho,0)}(\rho_1,z_1) = (\rho_{ra},z).$$
By the third property,
$$(\rho_1,z_1) = i_{(\varrho,0)} (\rho_{ra},z),$$
This means the set of all $(\rho_1,z_1)$ in \eqref{eq:InvComp} is the image of the line $\{(\rho_{ra},z)| z\in \mathbb{R}\}$
under the inversion $i_{(\varrho,0)}$, which is a circle. By symmetry,
this circle has a diameter on the $x$-axis. One end of the diameter is
$i_{(\varrho,0)} ((\rho_{ra}, \infty)) = (\varrho,0)$, and the other
end is $i_{(\varrho,0)} ((\varrho_{ra}, 0)) = \frac{R^2-1}{\varrho}$.
The circle is $\mathcal{C}(\varrho)$.
\qed

\gap
In virtue of Theorem~\ref{thm:FirstHomothetic}, we use the shorthand notation
$$
i_{\varrho}(T_R)
$$
to represent the common Euclidean shape of the cyclides in
\eqref{eq:SameShape}. Formally, $i_{\varrho}(T_R)$ is an element in the shape space \eqref{eq:ShapeSpace}.

\gap
To further analyze the shape $i_{\varrho}(T_R)$,
by Lemma~\ref{lemma:Maxwell} and Theorem~\ref{thm:FirstHomothetic}, it suffices to analyze the ratio
 $r_1:r_2:d$
of the cross-section of $i_{(\varrho,0,0)}(T_R)$
at its $P_1$ symmetry plane.

\begin{lemma} \label{lemma:CyclideInversions}
For any $R \in (1,\infty)$,
   the $P_1$ cross section of $\mathfrak{C} = i_{(\varrho,0,0)}(T_R)$ 
has the following measurements:
  \begin{itemize}
  \item[(i)] when $\varrho \in [0,R-1)$ (corresponding to the red circles in Figure~\ref{fig:Circles}(a)), the $P_1$ symmetry plane
of $\mathfrak{C}$ is the $x$-$z$ plane, and
  \bea
  \label{eq:ratioP1}
  \begin{aligned}
r_1:r_2:d = \lambda:1:\sqrt{(\lambda-1)^2 + 4 \lambda R^2}, \;\;
\mbox{where } \lambda = \frac{r_1}{r_2} = \frac{(\varrho+R)^2-1}{(\varrho-R)^2-1} \in [1,\infty).
   \end{aligned}
   \eea
  \item[(ii)] when $\varrho \in \left(R-1, \sqrt{R^2-1} \right]$ (corresponding to the blue circles in Figure~\ref{fig:Circles}(a)),
the $P_1$ symmetry plane
of $\mathfrak{C}$ is the $x$-$y$ plane, and
   \bea
  \label{eq:ratioP2}
  \begin{aligned}
r_1:r_2:d = \lambda:1:\sqrt{(\lambda-1)^2 + 4 \lambda \frac{R^2}{R^2-1}}, \;
\mbox{where } \lambda = \frac{r_1}{r_2} = \frac{(R-1)[(R+1)^2-\varrho^2]}{(R+1)[\varrho^2-(R-1)^2]}\in [1,\infty).
\end{aligned}
   \eea
\end{itemize}
\end{lemma}
\pf The first two steps of the proof of Theorem~\ref{thm:FirstHomothetic}
imply that
$$
P, \mbox{the $x$-$z$ plane,} \mbox{ is } \left\{
                               \begin{array}{ll}
                                 \mbox{the $P_1$ symmetry plane of $\mathfrak{C}$ when $\varrho \in [0,R-1)$} \\
                                 \mbox{the $P_2$ symmetry plane of $\mathfrak{C}$ when $\varrho \in \left(R-1, \sqrt{R^2-1} \right]$}
                               \end{array}
                             \right..
$$
In the first case, $\mathfrak{r}_i$ and $O_i$ in \eqref{eq:cyclideCrossSections} are such that
$\mathfrak{r}_1 > \mathfrak{r}_2$ and $O_2 < O_1$, and they give the $(r_1,r_2,d)$ measurements of $\mathfrak{C}$:
\bea \label{eq:measurements_i}
r_1  = \frac{1}{(\varrho-R)^2 -1}, \;\;
r_2 = \frac{1}{(\varrho+R)^2 -1}, \;\;
d = \frac{\varrho + R}{(\varrho + R)^2 -1} - \frac{\varrho - R}{(\varrho - R)^2 -1}.
\eea
In the second case, we also have $\mathfrak{r}_1 > \mathfrak{r}_2$ but now $O_1 < O_2$, and they give the
$(\tilde{r}_1,\tilde{r}_2,\tilde{d})$
measurements of $\mathfrak{C}$:
$$
\tilde{r}_1  = \frac{1}{1-(\varrho-R)^2}, \;\;
\tilde{r}_2 = \frac{1}{(\varrho+R)^2 -1}, \;\;
\tilde{d} = -\frac{\varrho + R}{(\varrho + R)^2 -1} + \frac{\varrho - R}{(\varrho - R)^2 -1}.
$$
By \eqref{eq:P1toP2}, we can convert the $(\tilde{r}_1,\tilde{r}_2,\tilde{d})$ measurements to
the $(r_1,r_2,d)$ measurements via $r_1 = (\tilde{r}_1-\tilde{r}_2+\tilde{d})/2$, $r_2 = (\tilde{r}_1-\tilde{r}_2-\tilde{d})/2$,
$d = \tilde{r}_1+\tilde{r}_2$, so
\bea \label{eq:measurements_ii}
r_1 = \frac{R-1}{\varrho^2-(R-1)^2}, \;\;
r_2 = \frac{R+1}{(R+1)^2-\varrho^2}, \;\;
d = \frac{1}{(R+\varrho)^2-1} - \frac{1}{(R-\varrho)^2-1}.
\eea
By routine computations, \eqref{eq:ratioP1} follows from \eqref{eq:measurements_i} and \eqref{eq:ratioP2} follows from \eqref{eq:measurements_ii}.
\qed

\gap
Lemma~\ref{lemma:CyclideInversions} has an almost immediate consequence:
\begin{theorem} \label{thm:MainStepI}
For any $R \in (1,\infty)$, $i_\varrho(T_R)$ is distinct for each $\varrho \in [0, R-1)$.
\begin{itemize}
\item If $R \neq \sqrt{2}$, then $i_\varrho(T_R)$ is distinct for each $\varrho \in [0, R-1) \cup (R-1, \sqrt{R^2-1}]$.
\item If $R = \sqrt{2}$, then 
$\varrho> \sqrt{2}-1$ adds no new shape and hence
the shape space \eqref{eq:ShapeSpace_Clifford} is in one-to-one correspondence
with
$$
 \left\{ i_{\varrho}(T_{\sqrt{2}}): \varrho \in [0, \sqrt{2}-1) \right\}.
$$
\end{itemize}
\end{theorem}
\pf Recall the two expressions in Lemma~\ref{lemma:CyclideInversions} for $\lambda = r_1/r_2$ in the two intervals of $\varrho$.
It is easy to check that both
$$\lambda_1: [0,R-1) \goto [1,\infty), \;\; \lambda_1(\varrho) = \frac{(\varrho+R)^2-1}{(\varrho-R)^2-1}$$
and
$$\lambda_2: \big( R-1,\sqrt{R^2-1} \big] \goto [1,\infty), \;\; \lambda_2(\varrho) =
\frac{(R-1)[(R+1)^2-\varrho^2]}{(R+1)[\varrho^2-(R-1)^2]}$$
are bijections: simply check that $\lambda_1$ is monotonic increasing from $1$ to $\infty$, and
$\lambda_2$ is monotonic decreasing from $\infty$ to $1$.
As the $r_1:r_2$ ratio of $i_\varrho(T_R)$ is distinct
for different $\varrho \in [0, R-1)$, the first statement of the theorem is true.
Likewise, $i_\varrho(T_R)$ is  also distinct for each $\varrho \in (R-1, \sqrt{R^2-1}]$.

To show the statement in the first bullet, it remains to argue that for $\varrho_1 \in [0, R-1)$ and $\varrho_2 \in (R-1, \sqrt{R^2-1}]$,
$i_{\varrho_1}(T_R) \neq i_{\varrho_2}(T_R)$. There are two cases:
\begin{enumerate}
  \item If $\lambda_1(\varrho_1) \neq \lambda_2(\varrho_2)$, then $i_{\varrho_1}(T_R) \neq i_{\varrho_2}(T_R)$.
  \item If $\lambda_1(\varrho_1) = \lambda_2(\varrho_2)$, then, by the expressions of the
$r_2:d$ ratio in Lemma~\ref{lemma:CyclideInversions},
the $r_2:d$ ratios of $i_{\varrho_1}(T_R)$ and $i_{\varrho_2}(T_R)$
are different exactly when $R^2 \neq \frac{R^2}{R^2-1}$. But
$$
R^2 = \frac{R^2}{R^2-1} \Longleftrightarrow R = \sqrt{2}.
$$
So we also have $i_{\varrho_1}(T_R) \neq i_{\varrho_2}(T_R)$ in this case.
\end{enumerate}
This argument proves the statement under the second bullet as well.
\qed

\gap
The next two results characterize the bigger shape space \eqref{eq:ShapeSpace}; they are inspiring for us but technically we do not need them
for this article. We
omit the detailed proofs, which
follow the same line of arguments as in that of Theorem~\ref{thm:MainStepI}.
\begin{lemma} \label{lemma:Duality}
For any $R \in (1,\infty)$, $\varrho \in [0, \sqrt{R^2-1}]$,
\bea
\label{eq:rhoDuality}
i_{\varrho}(T_{R}) = i_{\varrho'}(T_{R'}) \;\;\; \mbox{ where } \;\;\;
(R', \varrho') = \frac{1}{\sqrt{R^2-1}} \Big(R, \frac{\sqrt{R^2-1} - \varrho}{\sqrt{R^2-1}  + \varrho} \Big).
\eea
\end{lemma}

\begin{theorem}
Let
$$
C_R := \left\{
        \begin{array}{ll}
          \big[ 0, \sqrt{R^2-1} \big] \backslash \{R-1\} & \hbox{if $R \in (1,\sqrt{2})$} \\
          \big[ 0, \sqrt{2}-1 \big) & \hbox{if $R=\sqrt{2}$}
        \end{array}
      \right.,
\quad
C := \bigcup_{R \in (1,\sqrt{2}]} \{(R,\varrho): \varrho\in C_R\}.
$$
Distinct elements in $C$ 
correspond to distinct $i_{\varrho}(T_R)$ and
the shape space \eqref{eq:ShapeSpace} is in one-to-one correspondence
with
$$
\Big\{ i_{\varrho}(T_R): (R,\varrho) \in C \Big\}. 
$$
\end{theorem}

\begin{remark}
Admittedly, our proof of Theorem~\ref{thm:MainStepI} is very elementary
given the extensive development in M\"obius geometry; see, for example, \cite{Udo:MobiusBook,Jensen:Book}.
For instance, we use neither the usual representation of
the M\"obius group $\mbox{M{\"o}b}(3)$ in $\mathbb{S}^3$ nor its linear representation in the Lorentz space $\bR^{4,1}$.
It is unclear to us if our proof can be shortened
using the more modern techniques.
\end{remark}

\section{Step II: Rounding by sphere inversion} \label{sec:Step2}
\begin{theorem} \label{thm:Rounding}
If $S$ is a compact regular surface (with or without boundary) in $\bR^3$ and $p \in S$, then
\bea \label{eq:AreaRound}
{\rm Area}( i_q (S) ) \sim \frac{\pi}{|p-q|^2}, \;\; q \goto p, \;\; \overline{pq}\; \bot \; T_p S.
\eea
If $S$ is also closed and orientable (so $S$ and $i_q(S)$ have enclosing volumes), then
\bea \label{eq:VolumeRound}
{\rm Volume}( i_q (S) ) \sim \frac{\pi}{6|p-q|^3}, \;\; q \goto p, \;\; \overline{pq}\; \bot \; T_p S,
\eea \label{eq:IsoRound}
and (consequently)
\bea
v(i_q (S)) = \frac{{\rm Volume}( i_q (S) )}{(4\pi/3) \left({\rm Area}( i_q (S) )/(4\pi) \right)^{3/2}} \goto 1, \;\; q \goto p, \;\; \overline{pq}\; \bot \; T_p S.
\eea
\end{theorem}
\pf
Without loss of generality assume $p=(0,0,0)$ and $T_p S$ is the $x$-$y$ plane, and let $\varepsilon$ be a small scalar
representing the point $q=(0,0,\varepsilon)$ approaching the surface orthogonally at the origin. So the surface near $p$ can be written
as the graph of a smooth function $h(x,y)$, where $x^2+y^2<R^2$ for some $R>0$ and $h$ has a vanishing linear
approximation at the origin, i.e.
$h(0,0)=0 = \frac{\partial h}{\partial x}(0,0)=\frac{\partial h}{\partial y}(0,0)$,
and so
\bea \label{eq:h_prop}
h(x,y) = O(x^2+y^2), \;\; |\nabla h(x,y)| = O(\sqrt{x^2+y^2}), \quad (x,y)\goto (0,0).
\eea
Write $S_R := \{(x,y,h(x,y)): x^2+y^2<R^2\})$.
By continuity, the area of $i_{(0,0,\varepsilon)}(S \backslash S_R)$
approaches that of $i_{(0,0,0)}(S \backslash S_R)$ as $\varepsilon \goto 0$ and hence stays bounded for small $\varepsilon$.
So it suffices to prove \eqref{eq:AreaRound} with $S$ replaced by $S_R$.

The conformal factor of $i_{\mathbf{a}}$ is $\lambda^2(\mathbf{a},\mathbf{x}) = 1/\|\mathbf{x}-\mathbf{a}\|^4$, i.e.
$\langle d i_\mathbf{a}|_\mathbf{x} v, d i_\mathbf{a}|_\mathbf{x} w \rangle = \lambda^2(\mathbf{a},\mathbf{x}) \langle v,w \rangle$.
Therefore,
\begin{align} \label{eq:AreaInt}
\begin{split}
{\rm Area}( i_{(0,0,\varepsilon)} (S_R) )
& = \iint_{x^2+y^2<R^2} \frac{\sqrt{1+|\nabla h(x,y)|^2}}{[x^2+y^2+(h(x,y)-\varepsilon)^2]^2} \, dx dy \\
& = \int_{0}^{2\pi}
\left[
\int_{0}^{R} \frac{\sqrt{1+|\nabla h (r e^{i\theta})|^2}}{[r^2+(\varepsilon-h(r e^{i\theta}))^2]^2} \, r dr
\right]
d\theta.
\end{split}
\end{align}
Let $r_\ast(\varepsilon)=|\varepsilon|^{\alpha}$ for any $\alpha \in (1/2,1)$ so that
\bea \label{eq:rstar}
\mbox{(i)} \;|\varepsilon| = o(r_\ast(\varepsilon)) \;\;\mbox{ and }\;\; \mbox{(ii)} \; r_\ast(\varepsilon) =o( |\varepsilon|^{1/2}), \;\;\; \mbox{as $\varepsilon\goto 0$}.
\eea
We then split the inner integral in \eqref{eq:AreaInt} into
$\int_0^{r_\ast(\varepsilon)} + \int_{r_\ast(\varepsilon)}^R$; define
$$
J(\varepsilon) := \int_0^{2\pi} \int_0^{r_\ast(\varepsilon)} \frac{\sqrt{1+|\nabla h (r e^{i\theta})|^2}}{[r^2+(\varepsilon-h(r e^{i\theta}))^2]^2} \, r dr d\theta,
\quad
K(\varepsilon) := \int_0^{2\pi} \int_{r_\ast(\varepsilon)}^K \frac{\sqrt{1+|\nabla h (r e^{i\theta})|^2}}{[r^2+(\varepsilon-h(r e^{i\theta}))^2]^2} \, r dr d\theta.
$$
We shall prove \eqref{eq:AreaRound} by showing that the former integral is asymptotically equivalent to $\varepsilon^{-2}/2$ and the
latter grows slower than $\varepsilon^{-2}$.

For $J(\varepsilon)$, we compare it with the special case when $h\equiv 0$.
By \eqref{eq:h_prop}, there exists a constant $C>0$, independent
of $r$ and $\theta$, such that
$$|\nabla h(r e^{i\theta})|^2, \; |h(r e^{i\theta})| \leq C r^2.$$
For $r \in [0, r_\ast(\varepsilon)]$, $r^2\leq r_\ast(\varepsilon)^2 = o(|\varepsilon|)$ by \eqref{eq:rstar}(ii), so
$\varepsilon - h(re^{i\theta}) \sim \varepsilon$. Also, $1+|\nabla h (r e^{i\theta})|^2 \sim 1$.
From this it is easy to see that
\bea \label{eq:Jequiv}
J(\varepsilon)
\sim
\int_0^{2\pi} \int_0^{r_\ast(\varepsilon)} \frac{r}{[r^2+\varepsilon^2]^2}  \, dr d\theta.
\eea
The right-hand side is $J(\varepsilon)$ in the case of $h\equiv 0$, whose asymptotic can be easily determined:
\bea \label{eq:h0}
\int_0^{2\pi} \int_0^{r_\ast(\varepsilon)} \frac{r}{[r^2+\varepsilon^2]^2} \, dr d\theta =
\frac{2\pi}{\varepsilon^2} \int_0^{r_\ast(\varepsilon)/\varepsilon} \frac{s\, ds}{(1+s^2)^2}
=
\frac{2\pi}{\varepsilon^2} \left[ \frac{1}{2} - \frac{1}{2(1+(r_\ast(\varepsilon)/\varepsilon)^2)} \right]
\sim
\frac{\pi}{\varepsilon^2}, \;\; \varepsilon\goto 0.
\eea
In the last step above, we used \eqref{eq:rstar}(i).


For $K(\varepsilon)$, note that $\nabla h$ is bounded on a compact set, so
\begin{align*}
K(\varepsilon) =
\int_0^{2\pi} \int_{r_\ast(\varepsilon)}^R \frac{\sqrt{1+|\nabla h (r e^{i\theta})|^2}}{[r^2+(\varepsilon-h(r e^{i\theta}))^2]^2} \, r dr d\theta
& \leq 2\pi \int_{r_\ast(\varepsilon)}^R \frac{C}{[r^2]^2} \, r dr \\
& \leq 2\pi C \int_{r_\ast(\varepsilon)}^\infty r^{-3} dr = {\pi C} r_\ast(\varepsilon)^{-2} = o(\varepsilon^{-2}).
\end{align*}
In the last step above, we again used \eqref{eq:rstar}(i).

We have completed the proof of \eqref{eq:AreaRound}.

Let $B$ be a ball whose boundary is tangent to $S$ at $p$ and lies inside of $S$,
so ${\rm Volume} (B) \leq {\rm Volume} (S)$ and also
$${\rm Volume} (i_q (B)) \leq {\rm Volume} (i_q (S)).$$
As before, write $|p-q| = \varepsilon$.
Since $i_q (B)$ is a ball with diameter $\sim 1/\varepsilon$,
$$
{\rm Volume} (i_q (B)) \sim \frac{4\pi}{3}\left(\frac{1}{2\varepsilon}\right)^3 = \frac{\pi}{6\varepsilon^3}, \;\; \varepsilon \goto 0.
$$
So ${\rm Volume} (i_q (S))$ grows at least as fast as $\pi/(6\varepsilon^3)$. By the first part of the theorem and the isoperimetric inequality,
${\rm Volume} (i_q (S))$ cannot grow faster than $\pi/(6\varepsilon^3)$, and \eqref{eq:VolumeRound} is proved.
\eop

\gap
\begin{remark}
We thank I. Pinelis for the help in analyzing the asymptotic of the area integral \eqref{eq:AreaInt}; see
\url{https://mathoverflow.net/questions/353648/asymptotic-of-an-area-integral}. 
\end{remark}

\section{Step III: Reduction to P-recurrence} \label{sec:StepII}
In this section we express by P-recurrences the surface area and enclosing volume of $i_\mathbf{a}(T_{\sqrt{2}})$, where
$\mathbf{a}=[a,0,0]^T$, $a \in [0,\sqrt{2}-1)$, which are the same as those of
${\rm SCT}_\mathbf{a}(T_{\sqrt{2}})$. (Recall $i(T_{\sqrt{2}}) = T_{\sqrt{2}}$.)
From these, an associated P-recurrence related to the isoperimetric ratio of $i_\mathbf{a}(T_{\sqrt{2}})$
will also be derived.

\subsection{Area and volume integrals}
The conformal factor of a special conformal transformation ${\rm SCT}_\mathbf{a} := i \circ t_\mathbf{a} \circ i$ is
$$
\lambda^2(\mathbf{a},\mathbf{x}) = \frac{1}{(1+2\langle \mathbf{a}, \mathbf{x}\rangle + \langle \mathbf{a},\mathbf{a} \rangle \langle \mathbf{x}, \mathbf{x}\rangle)^2},
$$ i.e. $\langle dS_\mathbf{a}|_\mathbf{x} v, dS_\mathbf{a}|_\mathbf{x} w \rangle = \lambda^2(\mathbf{a},\mathbf{x}) \langle v,w \rangle$.
So the area and enclosing volume of ${\rm SCT}_{[a,0,0]}(T_{\sqrt{2}})$
are given by
\beq \label{eq:AreaVolume}
A(a) = \int_0^{2\pi} \int_0^{2\pi} Q(a; \mathbf{x})^{-2} \, d{\rm Area}(u,v),
\quad
V(a) = \int_0^1 \int_0^{2\pi} \int_0^{2\pi}  Q(a; \mathbf{x})^{-3} \, d{\rm Vol}(u,v,r),
\eeq
where
\beq
Q(a; \mathbf{x}) := \frac{1}{\lambda([a,0,0]^T,\bx)} =  1+2 \mathbf{x}_1 a  + \|\mathbf{x}\|^2 a^2,
\eeq
$$\mathbf{x}(u,v,r)=\left[ \big(\sqrt{2} + r\sin(v) \big)\cos(u), \;
\big( \sqrt{2} + r\sin(v) \big)\sin(u), \;
r\cos(v) \right], \;\;
u, v \in [0,2\pi], \;\; r \in [0,1],$$
$$
d{\rm Area}(u,v) = (\sqrt{2}+\sin(v)) \, du\, dv, \;\; d{\rm Vol}(u,v,r) = r(\sqrt{2} + r\sin(v)) \,du\, dv\, dr.
$$
Notice also that
$$\langle \mathbf{x}, \mathbf{x} \rangle = \| \mathbf{x} \|^2 = 2+ r^2 + 2\sqrt{2} r \sin(v).$$

\subsection{Holomorphic extension}
The integral definitions of $A$ and $V$ above extend from the interval $[0, \sqrt{2}-1)$ to a holomorphic function
on the open disk
$$
D := \{z \in \bC: |z| < \sqrt{2}-1\}. 
$$
To see this, note that the roots of $Q(z; \textbf{x})$, viewed as a quadratic polynomial in $z$, can be expressed as
$$
\frac{-\textbf{x}_1 \pm i \sqrt{\textbf{x}_2^2 + \textbf{x}_3^2}}{\|\textbf{x}\|^2},
$$
so their moduli are both $1/\|\textbf{x}\|$. But $\textbf{x}$ is a point in the boundary or interior of the solid torus $T$, so
$\|\textbf{x}\| \in [\sqrt{2}-1,\sqrt{2}+1]$, which is equivalent to $1/\|\textbf{x}\| \in [\sqrt{2}-1,\sqrt{2}+1]$.
This means
$$Q(z; \textbf{x}) \neq 0, \;\; \forall \, z \in D, \;\; \mathbf{x} \in T.$$
Therefore $Q(z; \textbf{x})^{-L}$, $L=2$ or $3$, is holomorphic in the first argument and continuous in the second.
A standard argument in complex analysis 
shows that
$A$ and $V$, defined based on the integrals in \eqref{eq:AreaVolume}, extend to holomorphic functions on $D$.

So from now on, we write $A(z)$ and $V(z)$ instead of $A(a)$ and $V(a)$.

\subsection{Power series at $z=0$}
Since $A$ and $V$ are even functions, the odd power Taylor coefficients at $z=0$ all vanish.
Denote by $a_j$ and $v_j$ the coefficients of $z^{2j}$ in the expansions of $A(z)$ and $V(z)$ at $z=0$, respectively.
An observation here is that
\begin{align} \label{eq:AV_TaylorCoeff}
\begin{split}
\frac{d^n A}{dz^n}(0) &= \iint_{\partial T} \frac{d^n}{dz^n} Q(z; \mathbf{x}(u,v,1))^{-2} \Big|_{z=0} \,d{\rm Area}(u,v) \\
\frac{d^n V}{dz^n}(0) &= \iiint_{T} \frac{d^n}{dz^n} Q(z; \mathbf{x}(u,v,r))^{-3} \Big|_{z=0} \,d{\rm Vol}(u,v,r),
\end{split}
\end{align}
and, {thanks to the evaluation at $z=0$}, the integrands above are \emph{polynomials} in $\textbf{x}_1$ and $\|\textbf{x}\|^2$, hence
are \emph{trigonometric polynomials} in $(u,v)$.

Using either \eqref{eq:AV_TaylorCoeff} or the generalized binomial theorem to
expand $Q(z; \mathbf{x})^{-L}$ into a power series of $z$, i.e.
$$
Q(z; \mathbf{x})^{-L} = \sum_{n=0}^\infty \binom{n+L-1}{n}
(-1)^n (2\mathbf{x}_1 z + \|\mathbf{x}\|^2 z^2)^n,
$$
together with the identity (of Wallis' integrals):
$$
\int_0^{2\pi} \cos^{n}(v)\, dv = \int_0^{2\pi} \sin^{n}(v)\, dv = \left\{
\begin{array}{ll}
\frac{2\pi}{2^n} \binom{n}{n/2}, & \hbox{$n$ even} \\
0, & \hbox{$n$ odd}
\end{array}
\right.,
$$
we have
\begin{align*}
\begin{split}
a_j &= \sum_{\ell=0}^j (-1)^{j-\ell}\, (j+\ell+1) \binom{j+\ell}{j-\ell} \,
 \int_0^{2\pi} \int_0^{2\pi} (2\mathbf{x}_1(u,v,1))^{2\ell} \|\mathbf{x}(u,v,1)\|^{2(j-\ell)} d{\rm Area}(u,v)
\\
&= \sum_{\ell=0}^j (-1)^{j-\ell}\, (j+\ell+1) \binom{j+\ell}{j-\ell} \, 4^\ell
 \underbrace{\int_0^{2\pi} \cos^{2\ell}(u)\, du}_{=2\pi\binom{2\ell}{\ell}/4^\ell} 
\underbrace{\int_0^{2\pi}(\sqrt{2}+\sin(v))^{2\ell+1} (3+2\sqrt{2}\sin(v))^{j-\ell} dv.}_{= \sum_{p=0}^{2\ell+1} \sum_{q=0}^{j-\ell}
\binom{2\ell+1}{p} \sqrt{2}^{2\ell+1-p} \binom{j-\ell}{q} 3^{j-\ell-q} (2\sqrt{2})^q \int_0^{2\pi} \sin^{p+q} (v) \,dv}
\end{split}
\end{align*}
So, \vspace{-.2in}
\begin{align}  \label{eq:A_Coeff}
\begin{split}
\quad a_j &= \sqrt{2}\pi^2 \sum_{\ell=0}^j (-1)^{j-\ell}\, (j+\ell+1) \binom{j+\ell}{j-\ell, \ell, \ell} \; \alpha_{\ell,j}, \\
\alpha_{\ell,j} &=  2^{\ell+2} 3^{j-\ell} \underbrace{\sum_{p=0}^{2\ell+1} \sum_{q=0}^{j-\ell}}_{p+q={\rm even}}
\binom{2\ell+1}{p} \binom{j-\ell}{q} \binom{p+q}{(p+q)/2} 2^{(q-3p)/2} 3^{-q}.
\end{split}
\end{align}
Similarly,
\begin{align*}
\begin{split}
v_j &= \sum_{\ell=0}^j (-1)^{j-\ell}\, \frac{(j+\ell+1)(j+\ell+2)}{2} \binom{j+\ell}{j-\ell} \,
 \int_0^1 \int_0^{2\pi} \int_0^{2\pi}(2\mathbf{x}_1(u,v,r))^{2\ell} \|\mathbf{x}(u,v,r)\|^{2(j-\ell)} d{\rm Vol}(u,v,r) \\
 &= \sum_{\ell=0}^j (-1)^{j-\ell}\, \frac{(j+\ell+1)(j+\ell+2)}{2} \binom{j+\ell}{j-\ell} \, 4^\ell
 \underbrace{\int_0^{2\pi} \cos^{2\ell}(u)\, du}_{=2\pi\binom{2\ell}{\ell}/4^\ell} \; \times \\
& \hspace{5cm} \underbrace{\int_0^1 \int_0^{2\pi} r (\sqrt{2}+r\sin(v))^{2\ell+1} (2+r^2+2\sqrt{2}r \sin(v))^{j-\ell} dv dr.}_{
\sum_{p=0}^{2\ell+1} \sum_{q=0}^{j-\ell}
\binom{2\ell+1}{p} \sqrt{2}^{2\ell+1-p}  \binom{j-\ell}{q} (2\sqrt{2})^q \int_0^{2\pi} \sin^{p+q} (v) \,dv \int_0^1 r^{p+q+1}(2+r^2)^{j-\ell-q}\, dr}
\end{split}
\end{align*}
So,  \vspace{-.2in}
\begin{align} \label{eq:V_Coeff}
\begin{split}
\quad v_j &=  \sqrt{2} \pi^2 \sum_{\ell=0}^j (-1)^{j-\ell}\, (j+\ell+1)(j+\ell+2) \binom{j+\ell}{j-\ell, \ell, \ell} \; \nu_{\ell,j}, \\
\nu_{\ell,j} &= 2^{\ell+1}  \underbrace{\sum_{p=0}^{2\ell+1} \sum_{q=0}^{j-\ell}}_{p+q={\rm even}}
\binom{2\ell+1}{p} \binom{j-\ell}{q} \binom{p+q}{(p+q)/2}\,2^{(q-3p)/2} \eta_{p,q,\ell,j}, \\
\eta_{p,q,\ell,j} &= \int_0^1 r^{p+q+1}(2+r^2)^{j-\ell-q}\, dr = \sum_{k}^{j-\ell-q} \binom{j-\ell-q}{k}
\frac{2^{j-\ell-q-k}}{2k+p+q+2}.
\end{split}
\end{align}
And we have the following power series:
\begin{align*}
\begin{split}
\frac{1}{\sqrt{2}\pi^2} A(z) =& \;
4 + 52\, z^2 + 477\, z^4 + 3809\, z^6 + \frac{451625}{16}\, z^8 +
\cdots\\
\frac{1}{\sqrt{2}\pi^2} V(z) = & \;
2 + 48\, z^2 + \frac{1269}{2}\, z^4 + 6600\, z^6 + \frac{1928025}{32}\, z^8 +
\cdots
\end{split}
\end{align*}
By the expressions \eqref{eq:A_Coeff}-\eqref{eq:V_Coeff},
$\frac{1}{\sqrt{2}\pi^2} a_n$, $\frac{1}{\sqrt{2}\pi^2} v_n$
are rational.

\subsection{Isoperimetric Ratio}
To show that the isoperimetric ratio of ${\rm SCT}_{[a,0,0]}(T_{\sqrt{2}})$
is monotonic increasing in $a \in [0,\sqrt{2}-1)$, it suffices to show
$$\Delta(a) := \frac{d}{da} \ln \frac{V(a)^2}{A(a)^3} = 2 \frac{V'(a)}{V(a)} - 3 \frac{A'(a)}{A(a)}>0,
\;\; \mbox{ or } \;\;
2 V'(a)A(a) - 3 V(a)A'(a)>0.
$$
It happens that $\Delta(a)$ is proportional to
the distance between the area and volume centers of the cyclide ${\rm SCT}_{[a,0,0]}(T)$. Precisely,
$\Delta(a) = 12 \big[ \mathbf{x}^A(a) - \mathbf{x}^V(a) \big]$
where $\mathbf{x}^A(a)$ and $\mathbf{x}^V(a)$ are the first coordinates of the area and volume centers
of ${\rm SCT}_{[a,0,0]}(T_{\sqrt{2}})$, respectively. This follows from
the observation that $({\rm SCT}_{[a,0,0]} \circ \mathbf{x})_1 = \frac{1}{2} Q'(a; \mathbf{x})/Q(a; \mathbf{x})$.

By the Taylor expansions of $A(a)$ and $V(a)$, we have
\begin{align} \label{eq:dk}
\begin{split}
&\frac{1}{2\pi^4} (2 V'(a) A(a) - 3 V(a) A'(a)) \\ 
=& \sum_k 
\Big[ \underbrace{2({v}_1 {a}_k + 2 {v}_2 {a}_{k-1} +
\cdots + (k+1)v_{k+1}a_0) -
3(a_1 v_k + 2 a_2 v_{k-1} + \cdots +
(k+1){a}_{k+1}{v}_0) }_{=: d_k} \Big]
a^{2k+1}\\
= &
72\, a + 1932\, a^3 + 31248\, a^5 + \frac{790101}{2}\, a^7 + \frac{17208645}{4} \, a^9+
\cdots
\end{split}
\end{align}

\subsection{P-recurrence} The combinatorial expressions
\eqref{eq:A_Coeff}-\eqref{eq:V_Coeff}, together with the closure properties of
holonomic sequences \cite{ZEILBERGER1990321,MR1676282,Kauers:book}, show
that $({a}_n)_{n\geq 0}$ and $({v}_n)_{n\geq 0}$ are
P-recursive,
i.e. they satisfy linear recurrences with polynomial coefficients.
Equivalently, their generating functions, namely
$$
\bar{A}(z) = \sum_{n\geq0} a_n z^n, \quad  \bar{V}(z) = \sum_{n\geq0} v_n z^n,
$$
are holonomic or $D$-finite, i.e. they satisfy linear differential equations with
polynomial coefficients. 
The generating functions of $({a}_n)_{n\geq 0}$ and $({v}_n)_{n\geq 0}$
are related to
the original area and volume functions $A(z)$ and $V(z)$ simply by
$A(z) = \bar{A}(z^2)$ and $V(z) = \bar{V}(z^2)$.
The generating function of the sequence $(d_k)_{k\geq 0}$, defined by \eqref{eq:dk}, is given by
\beq \label{eq:Dz}
\bar{D}(z):=
\sum_{n=0}^\infty d_n z^n = 2 \bar{V}'(z) \bar{A}(z) - 3 \bar{V}(z) \bar{A}'(z).
\eeq
Since holonomic functions are closed under Hadamard product (hence differentiation),
product, and linear combination, $(d_n)_{n\geq0}$ is also holonomic.

\begin{proposition}
The P-recurrences of $({a}_n)_{n\geq 0}$, $({v}_n)_{n\geq 0}$ and $(d_n)_{n\geq 0}$ are given by
\begin{align} \label{eq:a_recur}
\begin{split}
\sum_{i=0}^3 p_i(n) a_{n+i}=0, \mbox{ where }
\begin{bmatrix}
  p_0(n) \\
  p_1(n) \\
  p_2(n) \\
  p_3(n) \\
\end{bmatrix}
=
\begin{bmatrix*}[r]
-84 & -136 & -81 & -21 & -2\\
399 & 730 & 484 & 137 & 14\\
-474 & -835 & -529 & -143 & -14\\
54 & 99 & 66 & 19 & 2
\end{bmatrix*}
\begin{bmatrix}
  1 \\
  n \\
  n^2 \\
  n^3 \\
  n^4
\end{bmatrix}
\end{split}
\end{align}
\begin{align} \label{eq:v_recur}
\begin{split}
\sum_{i=0}^3 q_i(n) v_{n+i}=0, \mbox{ where }
\begin{bmatrix}
  q_0(n) \\
  q_1(n) \\
  q_2(n) \\
  q_3(n) \\
\end{bmatrix}
=
\begin{bmatrix*}[r]
-252 & -303 & -136 & -27 & -2\\
960 & 1384 & 730 & 167 & 14\\
-1008 & -1436 & -748 & -169 & -14\\
90 & 141 & 82 & 21 & 2
\end{bmatrix*}
\begin{bmatrix}
  1 \\
  n \\
  n^2 \\
  n^3 \\
  n^4
\end{bmatrix}
\end{split}
\end{align}
\begin{align} \label{eq:d_recur}
\begin{split}
&\sum_{i=0}^7 r_i(n) d_{n+i}=0, \mbox{ where } [r_0(n), r_1(n),\ldots,r_7(n)]^T = M [1, n, n^2, \ldots, n^7]^T, \\
M &=
\left[\begin{smallmatrix*}[r]
- \frac{1630207404}{1529} & - \frac{3176073675}{3058} & - \frac{660587685}{1529} & - \frac{1216898711}{12232} & - \frac{167529251}{12232} & - \frac{626799}{556} & - \frac{7141}{139} & -1\\
\frac{18219511026}{1529} & \frac{6798395835}{556} & \frac{16328931207}{3058} & \frac{15735207287}{12232} & \frac{2258693435}{12232} & \frac{8782801}{556} & \frac{103675}{139} & 15\\
- \frac{80949464718}{1529} & - \frac{338705850511}{6116} & - \frac{150907466733}{6116} & - \frac{74228837833}{12232} & - \frac{10882115811}{12232} & - \frac{43223443}{556} & - \frac{521157}{139} & -77\\
\frac{347623458975}{3058} & \frac{32991350565}{278} & \frac{322759355227}{6116} & \frac{158457515673}{12232} & \frac{23184921987}{12232} & \frac{91902509}{556} & \frac{1105723}{139} & 163\\
- \frac{368052969807}{3058} & - \frac{190572156372}{1529} & - \frac{168114763631}{3058} & - \frac{163720428321}{12232} & - \frac{23758375953}{12232} & - \frac{93404429}{556} & - \frac{1114663}{139} & -163\\
\frac{177327816597}{3058} & \frac{366011927673}{6116} & \frac{40230202855}{1529} & \frac{78121412337}{12232} & \frac{11304865929}{12232} & \frac{44328883}{556} & \frac{527737}{139} & 77\\
- \frac{29809040325}{3058} & - \frac{62775138251}{6116} & - \frac{28175845633}{6116} & - \frac{13970430847}{12232} & - \frac{2065443305}{12232} & - \frac{8275441}{556} & - \frac{100655}{139} & -15\\
\frac{818331696}{1529} & \frac{880217988}{1529} & \frac{1617383067}{6116} & \frac{822460415}{12232} & \frac{124982969}{12232} & \frac{515919}{556} & \frac{6481}{139} & 1
\end{smallmatrix*} \right].
\end{split}
\end{align}
Moreover, these are the only $P$-recurrences with the corresponding
order ($r$) and degree ($d$) for the three sequences. (E.g., \eqref{eq:a_recur} is
the only P-recurrence  with $(r,d)=(3,4)$ satisfied by the sequence $(a_n)$.)
\end{proposition}

A proof of the first part of the proposition, namely, the sequences defined by \eqref{eq:A_Coeff}-\eqref{eq:dk}
satisfy the P-recurrences
\eqref{eq:a_recur}-\eqref{eq:d_recur},
 can be established by a refinement of Zeilberger's \emph{creative telescoping method} \cite{ZEILBERGER1990321} due to
Koutschan \cite{Koutschan-2010} (implemented in his Mathematica package \verb$HolonomicFunctions$.)
Without diving into this method, we can check the second part of the claim in an elementary fashion. Assume that we have established that
$(a_n)$ follows a P-reccurence of order $r=3$ and degree $d=4$, then the $(d+1)(r+1)=20$ coefficients in the polynomials satisfy,
for every index $n$, a homogeneous linear equation with rational coefficients determined by the terms
$a_n, a_{n+1}, a_{n+2}, a_{n+3}$. Using the first $N+4$ terms of
the sequence $a_n$ with any $N\geq 20$, easily computable by \eqref{eq:A_Coeff}, we can set up a homogeneous linear system that must be satisfied by the 20 coefficients. Using a symbolic linear solver
to explicitly work out of a basis  of
the null space of the rational $N \times 20$ coefficient matrix -- and seeing that the basis
consists of one vector in $\bR^{20}$ with a certain $N\geq 20$ --
 would not only prove the claimed uniqueness (up to an arbitrary scaling factor), but also reproduce the
P-recurrence in \eqref{eq:a_recur}. This method is called `guessing' in \cite{Kauers:book}, as it
can be used to guess (with high confidence)
 what the P-recurrence might be when used with a big enough $N$.

Using asymptotic techniques \cite{Wimp-1985,Flajolet:2009:AC:1506267,Kauers:Mathematica} of
holonomic functions, it can also be shown that
\bea \label{eq:dn_asymp}
d_n \sim c \cdot \big(\sqrt{2}+1\big)^{2n} n^3 \ln(n), 
\quad
c \approx 8.071956... .
\eea
This is more than enough for showing that $d_n$ is eventually positive, but is insufficient for
verifying full positivity.

\section{Final Remarks} \label{sec:Final}
This paper connects a special case of the theory of Willmore surfaces to the theory of special functions, with the hope that it
may mobilize some interests in (i) the more ambitious uniqueness question discussed in
Section~\ref{sec:intro} and (ii) the positivity problem of P-recurrence, which is already
a well-known open problem in combinatorics \cite{Kauers_2010,Ouaknine_2013}. With all likelihood,
our approach for case (iii) of Conjecture~\ref{conjecture:Main}, being specific to the Clifford torus,
would not contribute much to the uniqueness problem in the other two cases Conjecture~\ref{conjecture:Main}. However,
it remains to see if the special function approach applies to the understanding of the
 higher genus Lawson surface $\xi_{g,1}$, which is conjectured to be the genus $g$ Willmore minimizer. (Recall that
the Clifford torus corresponds to $\xi_{1,1}$.)

On the other end, the ongoing work on attacking the positivity
of the P-recurrence \eqref{eq:d_recur} should contribute to the general positivity problem. A key difficulty of
proving positivity is that the characteristic polynomial of \eqref{eq:d_recur}, namely,
$z^7 - 15 z^6 + 77 z^5 - 163 z^4 + 163 z^3 - 77 z^2 + 15 z -1$ 
has roots
$$
\rho, \rho, 1, 1, 1, \rho^{-1}, \rho^{-1},
\quad \rho = \big(\sqrt{2}+1\big)^2.
$$
The repeated dominant root makes a certain dynamical system associated to the recurrence unstable,
which is related to why positivity is difficult to check.
In contrast, the characteristic polynomials of \eqref{eq:a_recur} and \eqref{eq:v_recur}, both being $z^3-7z^2+7z-1$, have roots
$\rho, 1, \rho^{-1}$; their positivity is easy to check by a simple inductive argument.

\bibliographystyle{plain}
  {\small
  \bibliography{refinement}
  }
\end{document}